\documentclass[10pt]{article}
\usepackage{amsfonts,amsmath,latexsym,amsthm,amssymb}
\usepackage{amstext}
\usepackage{xy,graphicx}
\input xy
\xyoption{all}
  
\newcommand{\HH}{\mathbb H}
\newcommand{\RR}{\mathbb R}
\newcommand{\ZZ}{\mathbb Z}
\newcommand{\NN}{\mathbb N}
\newcommand{\CC}{\mathbb C}
\newcommand{\DD}{\mathbb D}

\newcommand{\calA}{\mathcal A}
\newcommand{\calC}{\mathcal C}
\newcommand{\calH}{\mathcal H}
\newcommand{\calL}{\mathcal L}
\newcommand{\calO}{\mathcal O}

\newcommand{\calU}{\mathcal U}
\newcommand{\calV}{\mathcal V}
\newcommand{\calW}{\mathcal W}
\newcommand{\del}{\partial}
\newcommand{\olM}{\overline{M}}
\newcommand{\dxx}{\frac{dx}{x}}

\renewcommand{\sc}{{\mathrm{sc}}}
\newcommand{\fb}{{\mathrm{fb}}}
\newcommand{\fc}{{\mathrm{fc}}}
\newcommand{\II}{{\mathbb I}{\mathbb I}}
\newcommand{\calR}{\mathcal R}
\newcommand{\hB}{B}
  
\newcommand{\al}{\alpha}
\newcommand{\be}{\beta}
\newcommand{\de}{\delta}
\newcommand{\e}{\epsilon}

\newcommand{\frakp}{\mathfrak{p}}
\newcommand{\frakm}{\mathfrak{m}}
\newcommand{\frakd}{\mathfrak{d}}

\newcommand{\phg}{{\mathrm{phg}}}

\newcounter{fb}
\newcounter{fc}

\newtheorem{proposition}{Proposition}
\newtheorem{corollary}{Corollary}
\newtheorem{definition}{Definition}
\newenvironment{remark}{\paragraph{\it Remark.}}{\vskip0.4cm}

\begin{document}   
   
\title{Hodge cohomology of gravitational instantons}
\author{Tam\'as Hausel \\ UC Berkeley \and
Eugenie Hunsicker \\ Lawrence University \and
Rafe Mazzeo \thanks{Supported by the NSF grants
DMS-991975 and DMS-0204730}\\ Stanford University}
 
\maketitle   

\begin{abstract}
We study the space of $L^2$ harmonic forms on complete manifolds 
with metrics of fibred boundary or fibred cusp type. These metrics 
generalize the geometric structures at infinity of several different 
well-known classes of metrics, including asymptotically locally
Euclidean manifolds, the (known types of) gravitational instantons, 
and also Poincar\'e metrics on ${\mathbb Q}$-rank $1$ ends of locally 
symmetric spaces and on the complements of smooth divisors in K\"ahler 
manifolds. The answer in all cases is given in terms of intersection 
cohomology of a stratified compactification of the manifold. The $L^2$ 
signature formula implied by our result is closely related to the one 
proved by Dai \cite{dai} and more generally by Vaillant \cite{Va}, and 
identifies Dai's $\tau$ invariant directly in terms of intersection
cohomology of differing perversities. This work is also closely
related to a recent paper of Carron \cite{Car} and the forthcoming
paper of Cheeger and Dai \cite{CD}. We apply our results to a number
of examples, gravitational instantons among them, 
arising in predictions about $L^2$ harmonic forms in duality theories 
in string theory.

\end{abstract}

\section{Introduction}
The Hodge theorem for a compact Riemannian manifold $(M,g)$ identifies the
space $L^2\calH^*(M,g)$ of $L^2$  harmonic forms on $M$ with the de Rham
cohomology of this space. When $M$ is no longer compact, $L^2\calH^*(M,g)$
is 
still of considerable interest, but no general theorem identifies it
with a topologically defined group. In a number of special noncompact
geometric situations, there are topological interpretations of this
`Hodge cohomology' space.  These include the Hodge theorem for
manifolds with cylindrical ends in Atiyah-Patodi-Singer \cite{APS},
Cheeger's seminal work on Hodge theory on spaces with conic and iterated
conic singularities and its relationship with intersection theory
\cite{Ch1},
\cite{Ch2}, \cite{CGM}, the considerable literature on Hodge cohomology on
locally symmetric spaces, cf.\ in particular \cite{Z} and \cite{SS},
and the third author's work \cite{M1}, \cite{MPh} concerning
(asymptotically) geometrically finite hyperbolic quotients.

The aim of this paper is to prove a Hodge-type theorem for two different
classes of Riemannian manifolds, special cases of which arise frequently
in many interesting problems in geometry and mathematical physics.
These are fibred boundary and fibred cusp metrics.  Manifolds with fibred
boundary metrics include all identified classes of gravitational instantons,
the name coined by Hawking for complete hyperk\"ahler
four-manifolds.  Special cases of fibred cusp metrics include the
familiar `Poincar\'e' metrics in the theory of locally symmetric
spaces. Slightly more specifically, a product of a compact manifold with an
asymptotically locally Euclidean (ALE) manifold  is an example of a general
fibred boundary metric and a product of a compact manifold with a finite
volume hyperbolic cusp is an example of a fibred cusp metric, and
the most general case incorporates twisted versions of these examples, and 
also only requires the fibration structure to exist at the boundary. In
particular, there are two special and very familiar subclasses of metrics
amongst these: the ALE manifolds, also called scattering metrics,
and manifolds with asymptotically cylindrical ends, also called
$b$ metrics, which are fibred boundary and fibred cusp metrics,
respectively, with trivial fibre. We describe these rigorously
and in more detail below.


Let $\overline{M}$ be a smooth compact manifold with boundary, and suppose
that $x$ is a boundary defining function (thus $x$ vanishes on $\del
\overline{M}$ and $dx \neq 0$ there). We recall four classes of metrics in
terms of their behaviour in some neighbourhood $\calU$ of $\del
\overline{M}$. In the first two of these, $\overline{M}$ is arbitrary, but
in the latter two, we assume that $Y \equiv \del \overline{M}$ is the
total space of a fibration $\phi:Y \to B$ with fiber $F$.
\begin{itemize}
\item $g$ is called a $b$-metric on the interior $M$ of $\overline{M}$
if in $\calU$ it takes the form
\[
g = \frac{dx^2}{x^2} + h,
\]
where $h$ is a smooth metric on $\del \olM$ (i.e. nondegenerate up
to the boundary);

\item $g$ is called a fibred cusp metric if in $\calU$ it takes the
form
\[
g = \frac{dx^2}{x^2} + \tilde{h} + x^2 k,
\]
where $\tilde{h}$ is a smooth extension to $\calU$ of
$\phi^* h$, where $h$ is an arbitrary metric on $B$, and $k$
is a symmetric two-tensor on $\del M$ which restricts to a
metric on each fiber $F$;

\item $g$ is called a scattering metric if in $\calU$ it takes the
form
\[
g = \frac{dx^2}{x^4} + \frac{h}{x^2},
\]
where $h$ is a smooth metric on $\del \olM$;

\item $g$ is called a fibred boundary metric if in $\calU$ it takes
the form
\[
g =  \frac{dx^2}{x^4} + \frac{\tilde{h}}{x^2} + k,
\]
where $\tilde{h}$ and $k$ are as above.

\end{itemize}

We have made a simplification here in not allowing cross-terms in these
metrics, and members of these restricted classes are usually called exact
$b$-metrics, etc. This is not serious because as discussed
in the next section, Hodge cohomology is invariant under
quasi-isometries, and so these cross-terms can always be deformed
away without changing the Hodge cohomology. In fact,
we shall henceforth assume that a product structure $[0,1)_x \times \del Y$
is fixed on $\calU$ and that the metrics $h$ and $k$ in each of the four cases
are independent of $x$ with respect to this decomposition.
We shall simply write $h$ instead of $\tilde{h}$. This multi-warped
product structure simplifies computations and general
fibred boundary and fibred cusp metrics may be deformed to
ones of this type without affecting the Hodge cohomology.

These metrics, or at least special cases of them, are all familiar,
albeit in different coordinate systems. Thus if we set $x = e^{-t}$,
then a $b$-metric becomes $dt^2 + h$ on $\RR^+ \times \del M$, so
it has cylindrical ends, while the same change
of coordinates transforms a fibred cusp metric to
$dt^2 + h + e^{-2t}k$, which is a standard form for a
${\mathbb Q}$-rank $1$ cusp when $\del M$ is a torus bundle
over a torus. Similarly, if we set $x = 1/r$, then a
scattering metric becomes $dr^2 + r^2 h$ with $r \to \infty$,
which is the standard form of the infinite end of a cone,
and corresponds to the ALE class of gravitational instantons, such as the
Eguchi-Hanson metric. Finally, a fibred boundary metric transforms under
this coordinate change to $dr^2 + r^2 h + k$, which is a common form
for metrics in the ALF and ALG classes of gravitational instantons,
such as the Taub-NUT metric and reduced 2-monopole moduli space metric.

The obvious compactification of $M$ as the manifold with boundary 
$\overline{M}$ is useful for many purposes, but to state the Hodge 
theorems here we define a new compactification $X$ by collapsing the 
fibres $F$ of $\partial \olM$. When $F$ is a sphere, $X$ is a  
manifold, but in general $X$ is a stratified space with one singular 
stratum, which we denote $\hB$ (hopefully this should cause no
confusion), and principal stratum $M = X \setminus B$.
A neighbourhood of $B$ is a cone bundle with link $F$ over $B$.
In particular, when $B$ is trivial, $X$ is the one point compactification
of $M$, whereas when $F$ is trivial, $X=\olM$. In any case, we set
$b = \dim B$ and $f = \dim F$ throughout this paper. $X$ is called
a {\it Witt space} if $H^{f/2}(F) = 0$, and as we explain below, 
the analysis is much simpler in this case. 

Our main theorems relate the Hodge cohomology of $M$, with either a fibred 
boundary or fibred cusp metric, to the intersection cohomology of $X$. We 
refer to \S 2 for a review of these latter spaces and an explanation of the
notation in the following. 

\medskip

\noindent{\bf Theorem 1.}\ {\it Let $(M,g)$ be a manifold of dimension $n$ with
fibred boundary metric. Then for any degree $0 \leq k \leq n$,
there are natural isomorphisms
\[
L^2\calH^k(M,g) \longrightarrow
\left\{ 
\begin{array}{cll} 
& \mbox{\rm Im} \, \big(I\!H^k_{f+\frac{b+1}{2} -k}(X,\hB)
\longrightarrow I\!H^k_{f+\frac{b-1}{2} - k}(X,\hB) \big) \qquad
& \mbox{$b$ odd}  \\
&I\!H^k_{f+\frac{b}{2} -k}(X,\hB) & \mbox{$b$ even},
\end{array} \right.
\]
where the notation $I\!H_j^k(X, \hB)$ is explained in section 2.2.2, equation
(5).
}

\medskip

\noindent{\bf Theorem 2.}\ {\it Let $(M,g)$ be a manifold of dimension $n$ with
fibred cusp metric. Then for $0 \leq k \leq n$, there is a natural isomorphism
\[
L^2\calH^k(M,g) \longrightarrow \mbox{\rm Im} \, \left(
I\!H^k_{\underline{\frakm}}(X,\hB) \longrightarrow 
I\!H^k_{\overline{\frakm}}(X,\hB) \right)
\]
where $\underline{\frakm}$ and $\overline{\frakm}$ are the lower middle
and upper middle perversities. These give the same cohomology when $X$ is 
a Witt space, in which case we write simply 
\[
L^2\calH^k(M,g) \cong I\!H^k_{\frakm}(X,\hB).
\]
}

The perversity functions which arise in Theorem 1 are somewhat nonstandard, 
but they appear naturally in this problem. We shall return in another 
paper to a closer examination of the relationships between perversity 
functions and weighted $L^2$ cohomologies in these and other related 
geometric settings. However, for now note that an interesting special 
case occurs when $F$ is the sphere $S^f$, in which case $X$ is a manifold 
and  intersection cohomology reduces to ordinary cohomology. Then Theorem 1 
becomes

\begin{corollary}  Let $(M,g)$ be a manifold of dimension $n$ with a fibred
boundary metric where the fiber of $Y=\del M$ is a sphere; thus $M$
is identified with the complement of the submanifold $B$ in the
compact manifold $X$. Then for any degree 
$0 \leq k \leq n$, there are natural isomorphisms
\begin{eqnarray}
L^2\calH^k(M,g) 
\cong \left\{ \begin{array}{ll}
H^k(X, B) & k \leq \frac{b}{2} \\
H^k(X) & \frac{b}{2} <k < n-\frac{b}{2} \\
H^k(X \setminus B) & k\geq n-\frac{b}{2}
\end{array} \right.
\label{sphereeven}
\end{eqnarray}
if $b$ is even, and
\begin{eqnarray}
L^2\calH^*(M,g)
\cong \left\{ \begin{array}{ll}
H^k(X,B) & k \leq \frac{b-1}{2} \\
\mbox{\rm Im}\, \big(H^k(X,B) \longrightarrow H^k(X) \big)
\quad & k =
\frac{b-1}{2} +1 \\
H^k(X) & \frac{b+1}{2} <k < n-\frac{b+1}{2} \\
\mbox{\rm Im}\,\big(H^k(X) \longrightarrow H^k(X \setminus B)) & k = 
n-\frac{b+1}{2} \\
H^k(X \setminus B) & k\geq n-\frac{b-1}{2}
\end{array} \right.
\label{sphereodd}
\end{eqnarray}
if $b$ is odd.
\end{corollary}

The specialization of Theorem 2 is even simpler:
\begin{corollary}  Let $(M,g)$ be a manifold of dimension $n$ with a fibred
cusp metric where $F = S^f$ as in the previous corollary. Then
the compactification $X$ is a manifold and for any degree 
$0 \leq k \leq n$, 
\[
L^2\calH^k(M,g) = H^k(X).
\]
\end{corollary}

Two degenerate cases of Theorems 1 and 2 are fairly well-known:

\medskip

\noindent{\bf Theorem 1A.}\ 
{\it Let $(M,g)$ be a manifold of dimension $n$ with
scattering metric. Then there are natural isomorphisms
\[
L^2\calH^k(M,g) \longrightarrow
\left\{
\begin{array}{rlll}
& H^k(M,\del M) & \qquad & k < n/2, \\
& \mbox{\rm Im}\,\left(H^k(M,\del M) \to H^k(M)\right) & \qquad & k = n/2,
\\ 
& H^k(M) & \qquad & k > n/2.
\end{array}
\right.
\]
}

\medskip

\noindent {\bf Theorem 2A.}\ {\it Let $(M,g)$ be a manifold of dimension 
$n$ with $b$-metric. Then for any degree $0 \leq k \leq n$,
there is a natural isomorphism
\[
L^2\calH^k(M,g) \longrightarrow \mbox{\rm Im} \left( H^k(X,\hB) \to
H^k(X-\hB)\right)
\cong \mbox{\rm Im} \left( H^k(M,\del M)\to H^k(M)\right)
\]
}

Theorem 2A is proved in \cite{APS}, while Theorem 1A is stated 
in \cite{Me-scm}, but the proof does 
not seem to be readily available in the literature.  We prove these 
first as a warm-up to the more general cases because the proofs are 
structured similarly but present fewer analytic and geometric demands.

In all these results, but particularly in the latter two where
the notation is more familiar, it is apparent that the topological expressions
on the right depend on the stratification $(X,\hB)$, and not just on $X$.
The traditional hypotheses about perversities were designed to make the
corresponding intersection cohomology spaces independent of
stratification, but as explained in \S 2, this independence is lost
in certain degrees because of our use of slightly 
more general perversity functions. 

As already indicated, there is a simpler proof of Theorem 2 when $X$ 
is a Witt space. The reason is that with this hypothesis the range of
$d$ is closed in all degrees, and so the space of $L^2$ harmonic forms 
is isomorphic to the $L^2$ cohomology. One can then directly apply 
techniques of \cite{Z} which are mainly sheaf-theoretic and
topological. We discuss this further in \S 5.5.  Note that since
$L^2\calH^{n/2}(M,g)$ only depends on the conformal class
of $g$, we can also compute the middle degree Hodge cohomology
for fibred boundary metrics when $X$ is Witt. In fact, in
this case there is a trick to prove Theorem 1 in many cases:
if $k < n/2$ and $\omega \in L^2\calH^k(M,\hat{g})$ then $\omega \wedge \eta
\in L^2\calH^{n-k}(M \times S^{n-2k},\hat{g})$, in particular is a middle 
degree class, where $\hat{g}$ is the product of a fibred boundary metric on 
$M$ and the standard metric on the sphere, and where $\eta$ is the volume 
form on $S^{n-2k}$. The easier analytic argument now works provided 
the compactification $X \times S^{n-2k}$ is 
a Witt space, which requires that both $H^{(f+n)/2 - k}(F) = 
H^{k - (f+n)/2}(F) = 0$ (one of which is of course always true). 
This can always be
used, for example, to reduce Theorem 1A to a simple special case of 
Corollary 2 which follows easily from Theorem 2A. 
In the end, however, this would express the Hodge cohomology for a 
fibred boundary metric in terms of the homology of a 
different space altogether, hence is certainly less preferable. 

In any case, when $X$ is not a Witt space, one needs to do something
to confront the main issue that the range of $d$ is not closed.
The analytic machinery we introduce in \S4 and \S 5 provides one
avenue for doing this. Another possible approach involves
Carron's notion of non-parabolicity at infinity \cite{Car}.
In fact, Carron has used this method to characterize the 
$L^2$ cohomology of arbitrary complete Riemannian manifolds 
with flat ends. There is a substantial, but not complete, overlap of 
his results with ours; we comment on this further in \S 6.

Other work very closely related to ours is a forthcoming paper
by Cheeger and Dai \cite{CD} concerning the $L^2$ cohomology
of cone bundles. Since we have not yet seen this paper, we cannot
comment specifically on its relationship with the results here.
However, there seems to be substantial overlap; it is likely that
we could deduce some of their results using the methods here, and 
using parametrices in the edge calculus \cite{ma-edge}. 
Their methods should certainly give some of our results too.

Hodge theorems are of course closely related to index theorems,
and Theorems 1 and 2 imply a signature formula:
\begin{corollary} Let $(M,g)$ be a fibred boundary or fibred cusp metric. 
Then 
\[
\mbox{\rm{sgn}}_{L^2}(M,g) = {\rm{sgn}}\, \left(
\mbox{\rm Im}\,\big(
I\!H^*_{\underline{\frakm}}(X,B) \to 
I\!H^*_{\overline{\frakm}}(X,B)\big)\right).
\] 
\end{corollary}

This corollary is very closely related to the signature theorem for 
fibred cusp metrics proved by Dai \cite{dai} in a special case 
(using M\"uller's $L^2$ index theorem for manifolds with ends which are locally
symmetric of ${\mathbb Q}$ rank $1$), and in more generality by
Vaillant \cite{Va}. This theorem of Dai and Vaillant states that
\begin{equation}
\mbox{sgn}_{L^2}(M,g) = \mbox{sgn}\, \left(
\mbox{\rm Im}\,\big(H^*(M, \del M) \to H^*(M)\big)\right) + \tau,
\label{eq:daivaillant}
\end{equation}
where the final term is the $\tau$ invariant of the fibration
of $\del M$ defined by Dai \cite{dai}. Combining this with 
the above corollary gives the very interesting equality 
\begin{equation}
\tau = 
{\rm{sgn}}\, \left(\mbox{\rm Im}\,\big(
I\!H^*_{\underline{\frakm}}(X,B) \to 
I\!H^*_{\overline{\frakm}}(X,B)\big)\right)
-
\mbox{sgn}\, \left(
\mbox{\rm Im}\,\big(H^*(M, \del M) \to H^*(M)\big)\right).\label{eq:difftau}
\end{equation}
We discuss this further in \S 6 and \S 7, and shall explore this identity 
 in another paper. 

Our initial and primary motivation for this work came from predictions 
arising in duality theories in string theory, some of which we describe 
in \S 7. In particular, physicists have predicted the dimensions of
the spaces $L^2\calH^*$ on the moduli space of magnetic monopoles on 
$\RR^3$ \cite{sen1}, multi-Taub-NUT gravitational instantons \cite{sen2}, 
quiver varieties \cite{vafa-witten} and certain $G_2$ and 
$\mbox{Spin}(7)$ manifolds \cite{gomis-etal}. In many of these cases, the
metrics are of fibred boundary type, and our Theorem 1 confirms
most of these predictions.  A notable exception is the prediction 
for the $G_2$ manifold in \cite{gomis-etal}, which our results prove is 
false. Most of these results have been or could be proved by techniques 
available in the literature \cite{hitchin1} and \cite{Va}. In
particular, as we explain in \S 6 below, taking \cite{hitchin1}
into account, Dai's signature theorem \cite{dai} suffices to
calculate $L^2\calH^*$ for all hyperk\"ahler metrics of
fibred boundary type, i.e.\ all known gravitational instantons.
Our methods and results give a unified approach and has the advantage
of using only basic asymptotic properties of the metric, rather
than any refined properties, e.g.\ having a large symmetry group or 
special holonomy group.
Moreover, the interpretation of Hodge cohomology in terms of
the intersection cohomology of a compactification is 
very much in the spirit of the original Hodge theorem for compact manifolds. 
We hope \cite{compactifications} that the results here, as well as those in \cite{Hu}, 
suggest the correct form for a general result which would encompass
the remaining cases of these predictions. 

This paper is organized as follows. In \S 2 we review
$L^2$ cohomology and the basics of intersection theory,
focusing on spaces with only one singular stratum. We also define
two different versions of weighted $L^2$ cohomology. 
A review of the proof of the Hodge theorem for compact manifolds is
presented in the brief \S 3; this provides the basic analytic
structure for the proofs of our main theorems and we emphasize
here the main analytic points for which replacements are needed.
The Hodge theorems for $b$ and scattering metrics are proved
in \S 4; this is accompanied by a review of the requisite
analysis of $b$ pseudodifferential operators. The more general
Hodge theorems are proved in \S 5, first by identifying the Hodge 
cohomology with weighted cohomology, then by relating weighted cohomology 
to intersection cohomology.  We briefly explain the relationship of
our results to those of Dai, Vaillant, Cheeger, Hitchin and Carron 
in \S 6. Finally, in \S 7, we discuss the special cases of these
theorems which provided our original motivation, where
$M$ is one of the `gravitational instantons', of currency
in physics.

The authors wish to thank Jean-Paul Brasselet, Sergey Cherkis, Jaume Gomis, 
Nigel Hitchin, Richard Melrose and Andras Vasy for their interest and
advice. We are grateful also to both Xianzhe Dai and the referee for
reading the paper carefully and providing many good suggestions to improve
the exposition. E.H. thanks Stanford University for 
support. T.H. was supported by a Miller Research Fellowship at UC Berkeley.
R.M. was a visitor at MSRI, in the Spring of 2001, when this work was 
started. Finally, we would also like to acknowledge the hospitality of 
the Institute of Brewed Awakenings in Berkeley, where many of the ideas 
of this work percolated.

\section{Cohomologies}

We discuss various cohomology theories (in a loosely construed
sense) which play significant r\^oles in this paper.

As a general word about notation, if $\mathcal{F}$ is some function space
on the Riemannian manifold $(M,g)$, then $\mathcal{F}\Omega^*(M)$
denotes the space of sections of the exterior bundle $\bigwedge^*(M)$
with this regularity. When $\mathcal{F} = L^2$ or a Sobolev space, then we
indicate the dependence on the metric by writing $\mathcal{F}\Omega^*(M,g)$.

\subsection{$L^2$ and Hodge cohomology}
We start with a review of some facts about $L^2$ cohomology, and its
relationship to the space of $L^2$ harmonic forms.

The absolute cohomology $H^k(M)$ of a general (open) manifold $M$
is identified with the de Rham complex of smooth forms with unrestricted
growth at infinity,
\[
\ldots \longrightarrow \calC^\infty\Omega^{k-1}(M) \longrightarrow
\calC^\infty\Omega^k(M) \longrightarrow
\calC^\infty\Omega^{k+1}(M) \longrightarrow \ldots
\]
Similarly, its compactly supported cohomology $H^k_c(M)$ is computed by the
de Rham complex of smooth compactly supported forms,
\begin{equation}
\ldots \longrightarrow \calC^\infty_0\Omega^{k-1}(M) \longrightarrow
\calC^\infty_0\Omega^k(M) \longrightarrow
\calC^\infty_0\Omega^{k+1}(M) \longrightarrow \ldots
\label{eq:csc}
\end{equation}
It is well known \cite{dR} that these same cohomologies can also be
computed using the complexes of distributional forms $(\calC^{-\infty}
\Omega^*(M),d)$ and $(\calC^{-\infty}_0\Omega^*(M),d)$. However, there
are many interesting complexes incorporating restrictions on
regularity and growth at infinity between these extremes.
The most popular of these (for good reason) is $L^2$ cohomology in
the presence of a Riemannian metric. To define it, 
complete the differential complex (\ref{eq:csc}) with respect to the norms
on the exterior bundles associated to $g$ and the volume form $dV_g$ so as
to obtain the Hilbert complex
\begin{equation}
\ldots \longrightarrow L^2\Omega_g^{k-1}(M) \longrightarrow
L^2\Omega_g^{k}(M) \longrightarrow
L^2\Omega_g^{k+1}(M)\longrightarrow \ldots
\label{eq:l2c}
\end{equation}
Strictly speaking, this is not a complex since the differential $d$ is
defined at each stage only on a dense subspace. Thus the space of degree
$k$ should be defined as $\{\omega \in L^2\Omega^k(M,g): d\omega \in
L^2\Omega^{k+1}(M,g)\} \subseteq H^1\Omega^k(M,g)$. The cohomology of
(\ref{eq:l2c}) is called the $L^2$ cohomology of $M$, and denoted
$H^*_{(2)}(M,g)$. In other words,
\[
H^k_{(2)}(M,g) = 
\left\{\omega \in L^2\Omega^k(M,g): d\omega = 0 \right\}
\big/ \left\{d\eta: \eta \in L^2\Omega_g^{k-1}(M), \
d\eta \in L^2\Omega_g^{k}(M) \right\}.
\]

To set this into context, recall the Kodaira decomposition theorem, which
states that for arbitrary manifolds, there is an orthogonal decomposition
\begin{equation}
L^2\Omega^k(M,g) = L^2\calH^k(M,g) \ \oplus \
\overline{d\calC^\infty_0\Omega^{k-1}}
\ \oplus \  \overline{\delta\calC^\infty_0\Omega^{k+1}},
\label{eq:Kd}
\end{equation}
where the first summand consists of forms $\omega \in L^2\Omega^k(M,g)$
such that both $d\omega = \delta \omega = 0$. This is the space of $L^2$
harmonic fields, or Hodge cohomology, and is our main object of study.
The proof of the 
Kodaira decomposition is closely related to the essential self-adjointness
of $d+\delta$ on $L^2\Omega^*(M,g)$, which in turn follows from
Gaffney's $L^2$ Stokes theorem, cf.\ \cite{dR}. It follows from this that
the subspace of closed forms is precisely the sum of the first
two summands here, and hence
\begin{equation*}
\begin{split}
H^k_{(2)}(M,g) = & L^2\calH^k(M,g) \\
& \oplus \left\{d\eta \in L^2\Omega^k(M,g): \eta \in L^2\Omega_g^{k-1}(M)
\right\}/\overline{
\left\{d\eta \in L^2\Omega_g^{k}(M), \
\eta \in L^2\Omega_g^{k-1}(M)\right\}}.
\end{split}
\end{equation*}

In particular, when the range of $d$ from $L^2\Omega^{k-1}$ to
$L^2\Omega^k$ is not closed, then $H^k_{(2)}$ is infinite
dimensional. This behaviour occurs in many instances, e.g.\ on
Euclidean space, and indeed, is the reason for some of the
difficulties in understanding $L^2$ cohomology. However, we can
define the reduced $L^2$ cohomology
\[
\overline{H}^k_{(2)}(M,g) \ \ = \ \ \left\{\omega \in L^2\Omega^k(M,g):
d\omega = 0 \right\} \big/ \overline{\left\{d\eta \in L^2\Omega^{k}(M,g), \
\eta \in L^2\Omega^{k-1}(M,g) \right\}}.
\]
Combined with (\ref{eq:Kd}), this gives the useful isomorphism
\[
\overline{H}^k_{(2)}(M,g) \cong L^2\calH^k(M,g),
\]
which reveals the surprising fact -- certainly not apparent from
the basic definition -- that Hodge cohomology is invariant
under quasi-isometric changes of the metric. In other words, if
two metrics are comparable, $g' \leq c g$, $g \leq c' g'$, for
constants $c, c'>0$, then $\overline{H}^*_{(2)}(M,*)$ is the
same when computed with respect to either metric, and hence the same
is true for $L^2\calH^*(M,*)$. Moreover, if $H^k_{(2)}(M,g)$ is finite
dimensional, then it is naturally isomorphic to $L^2\calH^k(M,g)$.

Reduced $L^2$ cohomology is not quite as tractable as it might appear.
For example, it is quite important in calculations that there is a
Mayer-Vietoris sequence for unreduced $L^2$ cohomology,  
but this is true only in special cases for reduced $L^2$ cohomology.  

\subsection{Intersection cohomology}
We now review some definitions and facts about the intersection
cohomology of stratified spaces.

\subsubsection{Generalities}
Let $X$ be a stratified space of real dimension $n$ with no codimension
one singularities. We always assume, without further comment,
that this space satisfies some extra hypotheses: if a point $q\in X$ is
contained in the stratum of codimension $\ell$, then it has a neighbourhood
$\calU$ diffeomorphic to $\calV \times C(L)$, where $\calV$
is diffeomorphic to a Euclidean ball and is contained in that stratum
and $C(L)$ is the cone over a link $L$, which itself is a stratified
space (of dimension $\ell-1$).

A perversity $\frakp$ is an $n$-tuple of natural
numbers, $(p(1),p(2), \ldots, p(n))$ satisfying $p(1)=p(2)=0$ and
$p(\ell-1) \leq p(\ell) \leq p(\ell-1)+1$ for all $\ell \leq n$. Associated
to
such a space $X$ and perversity $\frakp$ is the intersection
complex $IC^*_\frakp(X)$, where, roughly speaking, the integer $p(\ell)$
regulates the dimension of the intersection of generic chains
with the stratum of codimension $\ell$. The homology of this complex
is the intersection homology $I\!H^\frakp_*(X)$. The dual intersection
cohomology $I\!H^*_\frakp(X)$ is more relevant to our purposes.

The following result is fundamental.

\begin{proposition}[\cite{GM2}] Let $X$ be a stratified space
and let $(\calL^*, d)$ be a complex of fine sheaves on $X$ with
cohomology $H^*(X,\calL)$. Suppose that if $\calU$ is a neighbourhood
in the principal (smooth) stratum of $X$, then $H^*(\calU,\calL) =
H^*(\calU,\CC)$, while if $q$ lies in a stratum of codimension
$\ell$, and $\calU = \calV \times C(L)$ as above, then
\begin{equation}
H^k(\calU,\calL) \cong I\!H_\frakp^k(\calU) =
\left\{ \begin{array}{ll}
I\!H_\frakp^k(L) & k \leq \ell-2-p(\ell) \\
0         & k \geq \ell-1-p(\ell)
\end{array} \right.
\label{eq:localcalc}
\end{equation}
Then there is a natural isomorphism between the hypercohomology
$\mathbb{H}^{\,*}(X,\mathcal{L}^*)$ associated to this complex
of sheaves and $I\!H^*_\frakp(X)$, the intersection cohomology of
perversity $\frakp$.
\label{pr:shchar}
\end{proposition}

Thus intersection cohomology with perversity $\frakp$ may be calculated 
using any fine sheaf, so long as its local cohomology satisfies 
(\ref{eq:localcalc}), which we refer to as {\it the local computation}.  
See also \cite{CGM} and \cite{Bo} for more on this.

This proposition is modified later in this section to provide
a link between weighted cohomology and intersection cohomology.

\subsubsection{Intersection cohomology for spaces with only two strata}
Suppose now that $X$ has only two strata: the principal smooth
stratum and the stratum of codimension $\ell$, which we denote $B$.
For convenience, we assume that $B$ 
is connected, although all results here generalize easily to allow $B$ 
to have many components 
(even of different dimensions, so long as their closures are
disjoint). We denote by $F$ the link associated to any point $q \in B$.
This is a smooth compact manifold of dimension $\ell-1$ with trivial
stratification, and $I\!H^*_\frakp(F)= H^*(F)$ no matter the perversity
$\frakp$. We associate to $X$ the manifold with boundary $\olM$ which
is obtained by blowing up $B$, i.e.\ replacing $B$ by
its spherical normal bundle. (This may be visualized as the complement
of a tubular neighbourhood of $B$ in $X$.) Notice that $\del M$
fibres over $B$ with fiber $F$.

The only part of the perversity which affects $IC^*_\frakp(X)$,
and hence $I\!H_\frakp^*(X)$, is the value $p(\ell)$. The basic
hypothesis on $\frakp$ implies that $0 \leq p(\ell) \leq \ell-2$,
and by (\ref{eq:localcalc}), only the spaces $H^k(F)$, $0 \leq k
\leq \ell-2-p(\ell)$, are relevant for the calculation of these
intersection spaces. We now introduce an extension of these definitions
by allowing $p(\ell)$ to take on any integer value. This does not
give anything dramatically new: when $p(\ell) \leq -1$,
the local calculations (\ref{eq:localcalc}) agree with those
for the computation of $H^*(X-B) = H^*(\olM)$, whereas when
$p(\ell) \geq \ell-1$ then they agree with those for the computation of
$H^*(X,B) \cong H^*(M,\del M)$. Thus for any $j \in {\mathbb Z}$
we fix the notation
\begin{equation}
I\!H^*_{j}(X,B) = \left\{ \begin{array}{ll}
H^*(X-B) & j \leq -1, \\
I\!H_\frakp^*(X) & 0 \leq j \leq \ell-2, \\
H^*(X, B) & j \geq \ell-1,
\end{array} \right.
\label{eq:extih}
\end{equation}
where in the middle case, $\frakp$ is {\it any} perversity with $p(\ell)=j$.

We note some properties of these extended groups. First, $I\!H^k_j(X,B)
\cong I\!H^{n-k}_{\ell-2-j}(X,B)$, just as with the standard intersection
cohomology groups. Next, suppose that $X$ is smooth and endowed with
the stratification $(X\setminus B,B)$,
where $B$ is just a distinguished smooth $(n-\ell)$-dimensional
submanifold. The link at any point $q \in B$ is $S^{\ell-1}$, and so if
$\calU = \calV \times C(S^{\ell-1})$ is a neighbourhood of a point $q
\in B$, then 
\[
I\!H_\frakp^k(\calU) =
\left\{\begin{array}{llc} &I\!H_\frakp^k(S^{\ell-1}) =
H^k(S^{\ell-1}) \qquad &k \leq \ell - 2 - p(\ell) \\
& 0  &k \geq \ell - 1 - p(\ell).
\end{array}\right.
\]
If $0 \leq p(\ell) \leq  \ell-2$, this equals $\CC$ for $k=0$
and $0$ for $k > 0$, which is the same local calculation as for
the ordinary cohomology of a smooth manifold; hence $I\!H_j^*(X,B) =
H^*(X)$ in this case. As expected, this is independent of the
submanifold $B$, hence of the choice of stratification of $X$,
because the perversity $\frakp$ is a `traditional' one.
However, in the other cases, when $j \leq -1$ or $j \geq \ell-1$,
$I\!H^*_j(X,B)$ depends strongly on $B$. We also remark that
this extension allows us to consider spaces with a codimension
one stratum, i.e.\ a boundary. In this case, the link of a point on the
boundary is any point, so the local calculations corresponding to
$j\leq -1$ and $j \geq 0$ give absolute and relative cohomologies,
respectively. 

This use of nonstandard perversities is now common in intersection
theory; for example, they enter into calculations of weighted cohomology 
on locally symmetric spaces \cite{Nair}.

Now return to the class of manifolds of interest in this paper
where $M$ is the interior of a compact $n$-dimensional manifold
with boundary $\olM$, such that $\del M$ is the total space
of a fibration, with base $B$ and fiber $F$,
$\dim B = b$, $\dim F = f$.  $M$ has two natural compactifications:
the first, $\olM$, is obtained simply by adding $\del M$,
while the second, $X$, is the result of collapsing the fibres
of $\del M$ in $\olM$. We write the image of $\del M$ in $X$ as $\hB$.
Thus $X$ is a stratified space with a single singular stratum, $\hB$,
of codimension $\ell = n-b$.

Let us calculate the extended intersection groups $I\!H^*_{j}(\olM,\hB)$.
The first step is to localize the calculation around $\hB$. Let
$N(\hB)$ denote a normal neighbourhood of $\hB$, so that
$X = M \sqcup N(\hB)$. The overlap $M \cap N(\hB)$ retracts onto
$\del M \cong \del N(\hB)$. For each $j$ there is a Mayer-Vietoris 
sequence 
\[
\longrightarrow I\!H^k_j(\olM,\hB) \longrightarrow H^k(M) \oplus
I\!H^k_j(N(\hB),\hB) \longrightarrow H^k(\del M) \longrightarrow.
\]
This is elementary since $M \cap N(B)$ retracts onto
a compact subset of $X \setminus B$. In any case, it suffices 
to calculate the groups $I\!H_j^k(N(\hB),\hB))$.

Assume $(\mathcal{L}_j^*, d)$ is a complex of fine sheaves, the
hypercohomology of which is isomorphic to $I\!H^*_j(N(\hB),\hB)$.
Choose an open cover $\{\calU_{\alpha}\}$ of $\hB$ in $X$ such that 
the bundle $\del M \to B$ is trivial over each $\calU_\alpha$;
this lifts to the cover $\underline{\calU}=\{\phi^{-1}
(\calU_{\alpha})\}$ of $N(\hB)$. The bigraded complex
of \v{C}ech cochains with coefficients in $\calL^*_j$
\[
\begin{xy}
\xymatrix{
\vdots &  \vdots &  \vdots & \\
C^0(\underline{\calU}, \mathcal{L}_j^2) \ar[u]^d \ar[r]^{\delta} &
C^1(\underline{\calU}, \mathcal{L}_j^2) \ar[u]^d \ar[r]^{\delta} &
C^2(\underline{\calU},\mathcal{L}_j^2) \ar[u]^d \ar[r]^{\delta} & \hdots\\
C^0(\underline{\calU}, \mathcal{L}_j^1) \ar[u]^d \ar[r]^{\delta} &
C^1(\underline{\calU}, \mathcal{L}_j^1) \ar[u]^d \ar[r]^{\delta} &
C^2(\underline{\calU},\mathcal{L}_j^1) \ar[u]^d \ar[r]^{\delta} & \hdots\\
C^0(\underline{\calU}, \mathcal{L}_j^0) \ar[u]^d \ar[r]^{\delta} &
C^1(\underline{\calU}, \mathcal{L}_j^0) \ar[u]^d \ar[r]^{\delta} &
C^2(\underline{\calU},\mathcal{L}_j^0) \ar[u]^d \ar[r]^{\delta} & \hdots
}
\end{xy}
\]
has hypercohomology which can be calculated using either of the two
associated spectral sequences, cf.\ \cite{BT}.  Consider first 
the spectral sequence which starts with 
with the vertical differential $d$. Any intersection of
neighbourhoods $\phi^{-1}(\calU_{\alpha})$ is of the form $(0,s)
\times F \times \calU'$, where $\calU'$ is an intersection of
the neighbourhoods $\calU_\alpha$ in $B$. The local calculation
(\ref{eq:localcalc}) gives that the $E_1$ term of the spectral sequence is:
$$
\begin{xy}
\xymatrix{
\vdots &  \vdots &  \vdots & \\
\quad \qquad \ 0 \quad \ar[r]^{\delta} \qquad  & \
\qquad  0 \quad \ar[r]^{\ \ \delta}
\quad & \quad \qquad 0 \quad \ar[r]^{\delta} \quad & \quad \ldots \\
C^0(\underline{\calU}, H^{\ell-2-j}(F)) \ar[r]^{\delta} &
C^1(\underline{\calU}, H^{\ell-2-j}(F)) \ar[r]^{\delta} &
C^2(\underline{\calU},H^{\ell-2-j}(F)) \ar[r]^{\ \ \delta} & \quad \hdots\\
\vdots &  \vdots &  \vdots & \vdots\\
C^0(\underline{\calU}, H^{1}(F)) \ar[r]^{\delta} &
C^1(\underline{\calU}, H^{1}(F)) \ar[r]^{\delta} &
C^2(\underline{\calU},H^{1}(F)) \ar[r]^{\ \ \delta} & \quad \hdots\\
C^0(\underline{\calU}, H^{0}(F)) \ar[r]^{\delta} &
C^1(\underline{\calU}, H^{0}(F)) \ar[r]^{\delta} &
C^2(\underline{\calU},H^{0}(F)) \ar[r]^{\ \ \delta} & \quad \hdots
}
\end{xy}
$$
In all the rows below level $\ell-1-j$, this is the same as the $E_1$ term
of the Leray-Serre spectral sequence for the bundle $\del M \to B$,
but all rows at level $\ell-1-j$ and above are set to zero.
The next differential, $d_1$, is the horizontal \v{C}ech differential
$\delta$. Using it to calculate the $E_2$ term gives a bigraded
diagram which agrees below level $j$ with the $E_2$ term of
the same Leray-Serre spectral sequence. The higher terms
$E_k$ of this truncated Leray-Serre spectral sequence
converge to the extended intersection cohomology $I\!H_j^*(N(\hB),\hB)$.

One can, for example, see by examining the further terms of the resulting
spectral sequence that this truncation does not change the
limit
$E^{p,q}$ for $p+q< \ell-1-j$.  Thus for $k < \ell-1-j$, 
$\sum_{p+q=k} E^{p,q} = I\!H^k_j(N(\hB),\hB) \cong H^k(\del M)$. Using
this in the  Mayer-Vietoris sequence, we find that for $k<\ell-1-j$,
$I\!H^k_j(\olM,\hB) \cong H^k(M)$.

\subsection{Weighted cohomology and intersection cohomology}
As we have already explained, one difficulty with $L^2$ cohomology 
is that in many cases the range of $d$ is not closed, and this leads 
to the (somehow spurious) infinite
dimensionality of the quotient spaces. There are many ways to circumvent
this, each of which involve a perturbation of the Hilbert
spaces $L^2\Omega^*(M)$. One possibility is to use an $L^p$ completion,
$p \neq 2$, cf.\ \cite{Z3}, but the loss of the Hilbert space structure
is unfortunate and unnecessary. An alternate and preferable
method involves the use of weighted $L^2$ norms. The associated
theory is called weighted cohomology.

We do not attempt a general definition of weighted cohomology, but
specialize directly to the cases of interest here. Thus let
$\olM$ be a compact  smooth manifold with boundary, with boundary
defining function $x$. If $a \in \RR$, then $x^a L^2(X)$ is the space
of all functions (or forms) $u = x^a v$ where $v \in L^2(X)$.
In the following, fix a fibred boundary metric $g_{\fb}$ and
a fibred cusp metric $g_{\fc}$ on $X$; we may as well assume that 
$g_{\fc} = x^2 g_{\fb}$. We also use the notation $\Omega^*_\fb$
and $\Omega^*_\fc$. This is explained in \S 5, but for now we
say only that this denotes the normalizations of the
exterior bundle corresponding to $g_\fb$ and $g_\fc$ which
are best suited for computations.

\begin{definition} For $a \in \RR$, define the Hilbert complexes
\begin{equation}
\ldots \to x^a L^2\Omega_\fc^{k-1}(M,g_{\fc}) \to x^a L^2\Omega_\fc^{k}
(M,g_{\fc}) \to x^a L^2\Omega_\fc^{k+1}(M,g_{\fc}) \to \ldots
\label{eq:cwa} 
\end{equation}
and
\begin{equation}
\ldots \to x^{a-1} L^2\Omega_\fb^{k-1}(M,g_{\fb}) \to
x^{a}L^2\Omega_\fb^{k}(M,g_{\fb}) \to
x^{a+1}L^2\Omega_\fb^{k+1}(M,g_{\fb})\to\ldots\label{eq:cwad}
\end{equation}
as completions of the de Rham complex of smooth compactly supported forms
with respect to the stated norms at each degree. We then set
$WH^k(M,g_{\fc},a)$ and $\calW H^k(M,g_{\fb},a)$ to be the cohomology of
these two complexes, respectively, {\bf at degree $k$}. Thus
\begin{equation}
WH^k(M,g_{\fc},a) = \frac{\left\{\omega \in x^a L^2\Omega_\fc^k(M,g_{\fc}):
d\omega = 0 \right\}}
{\left\{d\eta: \eta \in x^a L^2\Omega_\fc^{k-1}(M,g_{\fc}), \
d\eta \in x^a L^2\Omega_\fc^{k}(M,g_{\fc}) \right\}}
\label{eq:dwhk}
\end{equation}
and
\begin{equation}
\calW H^k(M,g_{\fb},a) = \frac{\left\{\omega \in x^{a} L^2
\Omega_\fb^k(M,g_{\fb}):
d\omega = 0 \right\}}{\left\{d\eta: \eta \in x^{a-1}
L^2\Omega_\fb^{k-1}(M,g_{\fb}), \
d\eta \in x^{a}L^2\Omega_\fb^{k}(M,g_{\fb}) \right\}}
\label{eq:dwhkp}
\end{equation}
\end{definition}

We will supress the metric in this notation when it is unambiguous.
Since fibred boundary and fibred cusp metrics are conformally
related, these two cohomologies are essentially the same.
More precisely,
\begin{equation}
\calW H^k(M,g_{\fb},a) = WH^k(M,g_{\fc},(n/2)-k+a).
\label{eq:wc1=wc2}
\end{equation}
Thus for the remainder of this section we discuss only $WH^k$.

Our main goal now is to relate the weighted cohomology for fibred
cusp metrics to intersection cohomology.

\begin{proposition} Suppose that $k-1+a-f/2 \neq 0$ for $0 \leq k \leq f$.
Then 
\[
WH^*(M, g_{\fc},a) \cong I\!H^*_{[a+ f/2]}(X, \hB),
\]
where $[a+f/2]$ is the greatest integer less than or equal to
$a+f/2$.
\label{pr:wcih}
\end{proposition}

\begin{proof}
We prove this by considering the complex of sheaves associated to 
$x^a L^2$, so that its hypercohomology equals $WH^*(M,a)$, and show 
that its entries satisfy the appropriate local 
calculation (\ref{eq:localcalc}). In order to apply Proposition 
\ref{pr:shchar}, however, we must first show that this sheaf is fine.

For each $k$, define the presheaf
\[
\calL^k_a(\calU) = \left\{
\begin{array}{lll}
&L^2\Omega^k(\calU) \qquad \qquad \quad
& \calU \cap \hB = \emptyset \\
& x^a L^2\Omega^k_\fc (\calU \setminus (\calU\cap \hB)) \quad &
\calU \cap \hB \neq \emptyset, \end{array}
\right.
\]
where the notation in this last line should be self-explanatory.
The associated sheaf is denoted $\calL_{a}^k$.

In general, the sheaf of (weighted) $L^2$ forms on the compactification
of a manifold is not fine unless one has a `good' partition of unity,
i.e.\ such that the cutoff functions $\chi_\alpha$ have gradients
which are bounded, uniformly in $\alpha$. However, such partitions
of unity are easy to construct for fibred cusp metrics, cf.\ the
essentially identical discussion in \cite{Z}. We construct
a good cover and partition of unity on $\olM$ as follows. First
choose a finite cover $\{\calU_\alpha\}$ of (the interior of)
$M$ such that all $j$-fold intersections of these sets are contractible.
Choose a similar cover $\{\calV_\beta\}$ of $B$, and let
$\calU'_\beta = \phi^{-1}(\calV_\beta) \times (0,\e)$, where
$\phi:\del M \to B$. Then for $\e$ sufficiently small,
$\{\calU_\alpha,\calU'_\beta\} = \{\calU''_\gamma\}$ is a good cover
for $\olM$. Now choose a partition of unity $\{\chi_\gamma''\}$
where the elements satisfy no additional extra requirements
over sets not intersecting $\olM$, but which have the form
$\phi^* \tilde{\chi}_\beta (y) \hat{\chi}(x)$ on neighbourhoods
intersecting the boundary; then it is easy to see that
$|d\chi''_\alpha| \leq C$ uniformly in $\alpha$, as required.

Now turn to the local cohomology computation, following the discussion 
in \S 2 of \cite{Z}. Over neighbourhoods not intersecting $\hB$ in $X$, 
we apply the standard Poincar\'e lemma. 
On the other hand, suppose that $\calU = \calV \times F \times
(0,\e)$, where $\calV \subset B$. First, by an adapted form
of the  K\"unneth theorem, the (weighted) $L^2$ cohomology of
the product neighbourhood $\calU$ is the same as the weighted
$L^2$ cohomology of $F \times (0,\e)$; this is valid since
the weight function $x^a$ does not depend on $y \in \calV$. This
reduces us to computing $WH^k(F \times (0,\e),dx^2/x^2 + x^2 k_F, a)$,
where for the moment we write the metric on $F$ as $k_F$ to distinguish
it from the form degree. Setting $r=-\log x$ to accord
with the notation of \cite{Z}, and regarding the weight
as a norm on the trivial local coefficient system, then the conclusion
of Corollary 2.34 from \cite{Z} is that:
\begin{itemize}
\item[i)] $WH^k(X,a)$ is finite dimensional, i.e.\ the denominator
in its definition is closed, if and only if the same is true for $WH^j((0,\e),
dx^2/x^2,k-j +a- (f/2))$ and simultaneously $H^{k-j}(F) \neq 0$, $j= 0,1$.
\item[ii)] If this condition is satisfied, then 
\[
WH^k(X,a) \ \cong 
\begin{array}{rl} & WH^0\big((0,\e),dx^2/x^2,k+a-f/2) \otimes H^k(F)\big) \\
&\ \oplus 
WH^1\big((0,\e),dx^2/x^2,k-1+a-f/2) \otimes H^{k-1}(F)\big).
\end{array} 
\]
\end{itemize}
In fact, this follows once again from the K\"unneth theorem in \cite{Z}. 
We have
\[
WH^0((0,\e),dx^2/x^2,b) = \left\{\begin{array}{lll}
& \CC  \qquad & b<0 \\  & 0 & \mbox{otherwise},
\end{array}\right.
\]
whereas $WH^1((0,\e),dx^2/x^2,b) = 0$ if $b\neq 0$ and is infinite 
dimensional when $b=0$ the range of $d$ is not closed at weight $0$).

Returning to the local calculation, we deduce that
\[
\begin{split}
WH^k(\calU,a) = \big(WH^0((0,\e), dx^2/x^2, k+a-f/2) \otimes H^k(F)\big) \\
\oplus \ \big(WH^1((0,s), dx^2/x^2, k-1+a-f/2) \otimes H^{k-1}(F)\big).
\end{split}
\]
Since we are assuming that $k-1+a-f/2 \neq 0$ when $0 \leq k \leq f$,
we obtain finally that
\[
WH^k(\calU, a) \cong WH^0((0,s),dx^2/x^2,k+a-f/2)
\otimes H^k(F) \cong
\left\{
\begin{array}{ll}
H^k(F) & k<f/2-a \\
0 & k \geq f/2-a.
\end{array}
\right. 
\]
Since the codimension of $\hB$ is $f+1$, this is the same as the
local calculation for $I\!H^*_\frakp(X)$ when $p(f+1)=[a+(f/2)]$.
\end{proof}

A closer reading of this proof, which we leave to the reader,
gives the following:

\begin{corollary} When $a$ is sufficiently large, then
\[
WH^k(M,g_\fc,a) = \calW H^k(M,g_\fb,a) = H^k(M)
\]
and
\[
WH^k(M,g_\fc,-a) = \calW H^k(M,g_\fb,-a) = H^k(M,\del M)
\]
for every $k = 0, \ldots, n$. If $F = \emptyset$, then
these equalities are true for any $a > 0$. 
\label{cor:luxembourg}
\end{corollary}

\subsection{Representing intersection cohomology with conormal forms}
It will be quite useful later to be able to represent
classes in intersection cohomology with forms which have some better 
regularity, especially near $\hB$ 
(or, in the other compactification of $M$, near
$\del \olM$). On a manifold with boundary, a natural and useful
replacement for smoothness at the boundary is `conormality'.
This is closely associated with $b$-geometry, which is discussed
in \S 4.1 and we refer ahead to that section for the definition
of the space of $b$-vector fields $\calV_b(M)$. For now, we 
say less formally that $V \in \calV_b$ if it
is a smooth vector field on $M$ and is tangent to $\del \olM$.
Let $\gamma \in \RR$, and define the space of conormal functions
of order $\gamma$ by
\[
\calA^\gamma(M) = \left\{u: |V_1 \cdots V_\ell u| \leq C x^\gamma\
\forall \, \ell\ \mbox{and}\ V_j \in \calV_b\right\}.
\]
Clearly, any conormal function is $\calC^\infty$ in the interior of $M$,
and it has full tangential regularity at the boundary. 
This definition extends directly to sections of vector bundles.

We now define a complex of conormal forms.  As we discuss later, cf.\ 
\S 5.1, the operator $d$ acting on $\Omega^*_\fc(M)$ involves
differentiations with respect to elements of $\calV_b$, but also
involves the nonsmooth term $x^{-1}d_F$. Hence we set
\[
\calA^a\Omega^k_{\fc,0}(M) = \{\alpha \in \calA^a\Omega^k_\fc(M):
d\alpha \in \calA^a\Omega^{k+1}_\fc(M))\}
\]
so that $(\calA^a\Omega^*_{\fc,0}(M),d)$ is a complex. In essence,
forms in this complex have a decomposition $\eta = \eta_0 + \eta'$
where $\eta' \in \calA^{a+1}\Omega^*_\fc(M)$ and $\eta_0
\in \calA^{a}\Omega^*_\fc(M)$ is fiber-harmonic as defined in \S 5.

It is well-known, cf.\ \cite{Me-aps} Prop. 6.13, that the relative and 
absolute cohomology of a manifold with boundary can be calculated 
using complexes of conormal forms. Generalizing this, we have
\begin{proposition}  The cohomology of the conormal complex
$(\calA^{a-f/2}\Omega_{\fc,0}^*(M),d)$ is isomorphic to the weighted
cohomology $WH^*(M,g_\fc,a)$. Provided $k-1+a -f/2 \neq 0$ for $0 \leq k 
\leq f$, it is also isomorphic to $I\!H^*_{[a+f/2]}(X)$.
\label{pr:ccwcic}
\end{proposition}

The argument to prove this is nearly identical to that for
Proposition \ref{pr:wcih}. The point is simply that $\calA^a
\Omega_{\fc,0}^*(M)$ is the space of global sections of a free 
sheaf, the local cohomology of which satisfies the same
local calculation as the sheaf of appropriately weighted $L^2$ forms.
We omit the details.

\section{Review of the compact Hodge theorem}

Despite some trade-off in work, we shall mainly use the
Hodge-deRham operator $D_g = d + \delta_g$, rather than its
square, $D_g^2 = \Delta_g$, the Hodge Laplacian. We now review one proof
of the standard Hodge theorem on compact manifolds which is phrased
in terms of $D_g$; this is intended as a guide for the analogous
arguments in the various noncompact settings considered below, and is
also meant to draw attention to certain analytic aspects of the argument
which are standard when $M$ is compact, but not so straightforward in these
other settings.

Recall the two most important components of the argument when $M$ is
compact. First, the ellipticity of the self-adjoint operator $D=d+\delta$
(we drop the subscript $g$ from now on) shows that it has a generalized
inverse $G$, which is a pseudodifferential operator of order $-1$. Write
$L^2\calH^*(M) = \mbox{ker}(D)$ and let $\Pi$ denote the orthogonal
projection $L^2\Omega^*(M) \to L^2\calH^*(M)$, so that $GD = DG = 
I-\Pi$. Implicit in this equation is the fact that the
kernel and cokernel of $D$ are both identified with $L^2\calH^*(M)$. We have
\begin{equation}
G: H^s \Omega^*(M) \longrightarrow H^{s+1} \Omega^*(M),
\qquad \Pi: H^s\Omega^*(M) \longrightarrow \calC^\infty\Omega^*(M),
\qquad \mbox{for all}\ s \in \RR,
\label{eq:fredcpct}
\end{equation}
and of course $\Pi$ is finite rank. Also, both $d$ and $\delta$
both commute with $G$. It follows directly that $D$ is
Fredholm, e.g.\ on $L^2\Omega^*(M)$. Furthermore, (\ref{eq:fredcpct})
also shows that the deRham cohomology $H^k(M)$ can be calculated
using any one of the complexes, $\calC^{\infty}\Omega^*(M)$,
$L^2\Omega^*(M)$, or $\calC^{-\infty}\Omega^*(M)$, i.e., of smooth,
$L^2$ or distributional (current) forms.

Now to the argument. First let $\omega \in L^2\calH^k(M)$. Then $D^2\omega =
0$,
and since $\omega$ is smooth, there is no problem in the integration by
parts, 
\[
\langle D^2 \omega,\omega \rangle = \langle D\omega, D\omega \rangle =
\|d\omega\|^2 + \|\delta \omega\|^2,
\]
so that $d\omega = \delta \omega = 0$. In particular, $\omega$ is closed
and $[\omega]\in H^k(M)$ is well-defined. This defines a map
\[
\Phi: \calH^k(M) \longrightarrow H^k(M).
\]
We must show that $\Phi$ is both injective and surjective.

Suppose $\Phi(\omega) = [\omega] = 0$, i.e.\ $\omega = d\zeta$
for some $(k-1)$-form $\zeta$. We may assume that we are calculating
using the complex of smooth forms, so we can choose $\zeta$ to be
smooth.  Since there are no boundary terms to worry about, we can
integrate by parts to obtain:
\begin{equation}
||\omega||^2 = \langle \omega, d\zeta \rangle = \langle
\delta \omega, \zeta \rangle = 0, \label{eq:injg}
\end{equation}
and so $\Phi$ is injective.

Next, let $[\eta] \in H^k(M)$ and choose a smooth representative
$\eta$. Applying $GD = I-\Pi$ to it yields
\begin{equation}
\eta = D \zeta + \gamma, \qquad \mbox{where}\
\zeta = G\eta,\ \gamma = \Pi \eta.
\label{eq:cok}
\end{equation}
By (\ref{eq:fredcpct}) again, $\zeta \in \calC^\infty\Omega^*$, and of
course the same is true for $\gamma \in L^2\calH^*$.
Because $D$ and $G$ act on forms of all degrees together, we do not
know yet that $\zeta$ or $\gamma$ are forms of pure degree $k-1$ and
$k$, so we argue as follows. Write
\[
\delta \zeta = \eta - d\zeta - \gamma;
\]
then
\[
\|\delta \zeta\|^2  = \langle \delta \zeta, \eta - d\zeta - \gamma \rangle
= \langle \zeta, d\eta - d^2 \zeta - d\gamma\rangle = 0.
\]
Hence $\eta = d\zeta + \gamma$ where $\gamma \in L^2\calH^*(M)$.
Now, clearly, neither $\zeta$ or $\gamma$ have terms of degree
other than $k-1$ and $k$, respectively.
This establishes surjectivity of $\Phi$ and completes the proof.

\medskip

When $(M,g)$ is noncompact, each of these steps may fail in a variety
of ways, and our main task is to show that they can be justified for
fibred boundary and fibred cusp metrics. Most fundamentally, $D$ may
no longer be Fredholm on $L^2\Omega^*$, and so we must find some other
function space on which it does have closed range. In fact, in our
cases it is Fredholm on a scale of weighted $L^2$ spaces, and we must
study the action of $D$ on these spaces. In particular, we wish to find
function spaces $X$ and $Y$ such that $D:X\to Y$ is Fredholm with
cokernel identified with $L^2\calH^*(M)$. This will serve as the
replacement for (\ref{eq:cok}). To justify the various integrations
by parts, we must also establish that elements of $L^2\calH^*(M)$ decay
at some definite rate at infinity, and also show similar decay and
regularity properties for the forms $\zeta$.

\section{Nonfibred ends}

Although the $L^2$ Hodge theorems for $b$ (cylindrical) and scattering 
(asymptotically Euclidean) metrics, Theorems 1A and 2A, are already known,
we nevertheless present proofs of these results here which
address some (but not all) of the difficulties encountered in the general
fibred boundary and fibred cusp cases. 

We shall sometimes denote the Hodge-de Rham operator for a $b$ or
scattering metric by $D_b$ or $D_\sc$, respectively.
Recall from the end of the last section that we need to find
function spaces on which these operators have closed range,
and we must also establish various decay and regularity properties
for the $L^2$ harmonic forms, as well as the other auxiliary
forms which enter into the proof. To obtain these properties,
we use the machinery of the $b$-calculus \cite{Me-aps}, cf.\
also \cite{ma-edge}. In other words, we adopt the point of view that
in either case $D$ is an elliptic element in an appropriate ring
of degenerate differential operators on the manifold $\olM$.
Mapping and regularity properties of these operators can be
investigated using a parametrix for $D$ constructed in an
associated calculus of degenerate pseudodifferential operators.

\subsection{$b$ metrics and operators}

Let $g$ be an exact $b$ metric. Associated to it is the space of $b$
vector fields $\calV_b$, which by definition is the Lie algebra of
all smooth vector fields on $\olM$ which are tangent to $\del M$. In
a coordinate chart $(x,y_1, \ldots, y_{n-1})$ near $\del M$, where
$(y_1,\ldots, y_{n-1})$ are coordinates on $\del M$ extended to the
collar neighbourhood $\calU$ and $x$ is a boundary defining function,
any $Z \in \calV_b$ can be written as
\[
Z = a(x,y)x\del_x + \sum_{j=1}^{n-1} b_j(x,y)\del_{y_j}, \qquad a, b_j \in
\calC^\infty(\olM).
\]
Notice that $\calV_b$ contains precisely those smooth vector fields on
$\olM$ which have pointwise bounded norms with respect to any $b$-metric.
The vector fields $x\del_x$, $\del_{y_j}$ form a local spanning set
of a vector bundle over $M$, called the $b$-tangent bundle,
${}^bTM$. This bundle is canonically isomorphic to the ordinary
tangent bundle $TM$ only over the interior, $M$, of $\olM$, but the
canonical map ${}^bTM \to TM$ given by evaluating sections at a
point extends to a map which is neither injective nor surjective over
$\del M$; its nullspace is one-dimensional and is spanned by $x\del_x$.
The dual of
${}^bTM$ is the
$b$-cotangent bundle, ${}^bT^*M$, which is locally spanned by the
one-forms $dx/x$,  $dy_j$.
We write ${}^b\bigwedge^*M$ and $\calC^\infty\Omega^*_b(M)$ for the
exterior powers of this bundle and its space of smooth sections,
respectively. 

A differential operator $P$ on $M$ is called a $b$-operator if it
can be written locally as a sum of products of elements of $\calV_b$.
Thus, in these coordinates,
\[
P = \sum_{j+|\al|\leq m} a_{j,\al}(x,y)(x\del_x)^j\del_y^\al,
\]
with all coefficients $a_{j,\al}\in\calC^\infty(\olM)$. If $P$ is an
operator
on a space of sections of a bundle over $M$, then the coefficients
$a_{j,\al}$ will be smooth endomorphisms of the bundle.  The
$b$-symbol 
\[
{}^b\sigma_m(P)(x,y;\xi,\eta) =
i^{-m}\sum_{j+|\al|=m}a_{j,\al}(x,y)\xi^j \eta^\al
\]
is invariantly defined as a homogeneous function on ${}^bT^*M$, and
$P$ is elliptic in this setting if ${}^b\sigma_m(P)$ is nonvanishing
(or invertible if $P$ is a system) for $(\xi,\eta) \neq 0$.

Our primary example of a $b$ differential operator is the Hodge-deRham
operator $D = d + \delta$ with respect to a $b$ metric $g$ on
$M$. To illustrate the definitions above, we determine its form now,
assuming that the metric $h$ which appears in the decomposition
of $g$ does not depend on $x$ in the boundary neighbourhood $\calU$.

Near $\del M$ any element of $\Omega^k_b(M)$ can be written as
\[
\omega = \al + \frac{dx}{x}\wedge \be,
\]
where $\al(x,y)$ and $\be(x,y)$ are families of $k$- and $(k-1)$--forms,
respectively, on $\del M$ depending on $x$ as a smooth parameter.
The $L^2$ norm is given by
\[
||\omega||^2 = \int_{M} \left(|\al|_h^2 + |\be|_h^2\right) \,
\frac{dx\,dy}{x}.
\]

Since $b$-metrics are special cases of fibred cusp metrics,
where the fibration $\del M \to B$ has trivial fibres,
we cohere with the more general notation of this paper
and identify $\del M$ with $B$. The induced differential
is written $d_B$, and the codifferential, induced by
the metric $h$ on $B$, is written $\de_B$. We have
\begin{align}
d\omega & = d_B \al + \dxx \wedge (x\del_x \al - d_B \be),
\label{eq:dob}  \\
\de \omega &= \de_B \al - x\del_x \be - \dxx \wedge \de_B \be.
\label{eq:dstob}  
\end{align}
Finally, the $b$ symbol of $D$ is computed just as in the standard case,
so that if $\zeta = (\xi,\eta) \in {}^bT^*M$, then ${}^b\sigma_1(D) =
i\, (\zeta \wedge \cdot + \, \iota_\zeta\cdot )$. This gives
\begin{proposition} The operator
\[
D = d + \de: \calC^\infty\Omega^*_b(M) \longrightarrow
\calC^\infty\Omega^*_b(M)
\]
on $(M,g)$ is an elliptic $b$-differential operator of order $1$.
\end{proposition}

\begin{remark} It is natural to write forms in terms of the $b$ covector 
fields $dy_j$ and $\dxx$, since these have (essentially) unit length, 
but it is also important, since in a poorly chosen coframe, the
expression of $D$ might no longer be a $b$ operator. For example,
this is the case if we use the standard basis $dx$ and $dy_j$.
\end{remark}

Unlike the usual interior calculus, symbol ellipticity alone is
not enough to determine whether a $b$ differential operator $P$ is
Fredholm. For this one must also use another model for $P$, called the
indicial operator $I_P$. This operator acts on functions on $B_y \times 
\RR^+_s$ and is invariant with respect to dilations in $s$; for a 
general $P$ written as above, 
\[
I_P = \sum_{j+|\al| \leq m} a_{j,\al}(y)(s\del_s)^j \del_y^\al.
\]
To analyze this operator we use its dilation invariance. Thus
conjugating $I_P$ by the Mellin transform in $s$, 
\[
u(s,y) \longmapsto u_M(\gamma,y) = \int_0^\infty s^{\gamma}
u(s,y)\,\frac{dsdy}{s}, \qquad \gamma \in \CC,
\]
yields the indicial family, $I_P(\gamma)$, which is a holomorphic family of
elliptic operators on $B$ (when $P$ is $b$-elliptic). By the analytic
Fredholm 
theorem, this family is either never invertible for any $\gamma$ or else is
invertible for all $\gamma \in \CC \setminus \Lambda$, where $\Lambda$ is
a discrete set of complex numbers, called the indicial set, the elements of
which are called the indicial roots of $P$. It is not hard to see that the
first possibility never holds.
We shall use an alternate (equivalent) characterization of this indicial set
which is more intuitive, and certainly easier to calculate:
\[
\gamma \in \Lambda \Longleftrightarrow\ \exists \, \phi \in
\calC^\infty(Y)\ \mbox{such that}\
P(x^\gamma \phi(y)) = {\mathcal O}(x^{\gamma+1}) \qquad
\mbox{where $Y = \del M$}.
\]
Notice that $P(x^\gamma \phi(y)) = {\mathcal O}(x^\gamma)$ for
all $\gamma$ and $\phi$, and so $\gamma \in \Lambda$ if and only if there
is some additional cancellation, which arises precisely when
there is an element $s^\gamma \phi(y)$ in the nullspace of $I_P$.

Again we illustrate this through the operator $D$. Since we
are assuming that $h$ does not depend on $x$ in $\calU$, we can
identify $I_D$ with $D$ near $\del M$, and so all approximate
solutions of $D\omega = 0$ in the sense above are exact solutions
in this boundary neighbourhood. Now write $\omega = \omega' x^\gamma$, where
$\omega' = \alpha' + \dxx \wedge \beta'$ and neither $\alpha'$
nor $\beta'$ depend on $x$. Then by (\ref{eq:dob}), (\ref{eq:dstob}),
in $\calU$, 
\begin{equation}
D(\omega' x^\gamma) \equiv x^\gamma I_D(\gamma)(\omega') =
x^\gamma \left(D_B \al' - \gamma \be' + \dxx \wedge ( \gamma \al' - D_B \be')\right).
\label{eq:b9.5}
\end{equation}
Hence $\gamma$ is an indicial root if and only if there is
a solution $\omega'$ of the equations
\begin{equation}
D_B \al' = \gamma \be', \ D_B \be' = \gamma \al',
\label{eq:b10}
\end{equation}
which implies
\begin{equation}
\Delta_B \al' = \gamma^2 \al',\ \Delta_B \be' = \gamma^2 \be'.
\label{eq:b11}
\end{equation}
Thus $\gamma$ is an indicial root of $D$ if and only if $\gamma^2 \in
\mbox{spec}\,(\Delta_B)$ on $\Omega^*(B)$. Note that the operators in
(\ref{eq:b11}) preserve the form degree, and so are easier to analyze
than the operators in (\ref{eq:b10}). However, arbitrary solutions of
(\ref{eq:b11}) do not necessarily satisfy (\ref{eq:b10}); in other words, 
we must be cautious not to introduce spurious indicial roots by 
all solutions of the decoupled equations. From the Kodaira decomposition on
$\Omega^*B$, the only coupling in (\ref{eq:b10}) is between closed
$k$-forms and coclosed $(k-1)$-forms for each $k$. Thus let $\phi_j$ and 
$\psi_j$ be a complete set of eigenforms for $\Delta_B$ on coclosed 
$(k-1)$- and closed $k$-forms, with 
eigenvalue $\lambda_j^2$ and such that $d_B \phi_j = \lambda_j \psi_j$,
$\delta_B \psi_j = \lambda_j \phi_j$ for $\lambda_j \neq 0$.
Writing
\[
\al' = \sum \al_j(x)\psi_j, \qquad \be' = \sum \be_j (x)\phi_j,
\]
then (\ref{eq:b10}) gives
\[
\gamma \al_j = \lambda_j \be_j, \quad \gamma \be_j = \lambda_j \al_j,
\] 
which implies $\gamma^2 = \lambda_j^2$, as expected. We see finally that
\begin{equation}
\omega' = \sum_j \left\{\left(\al^+_j + \dxx \wedge
\be^+_j\right)x^{\lambda_j}
+ \left(\al^-_j + \dxx \wedge \be^-_j\right)x^{-\lambda_j}\right\},
\label{eq:b13}
\end{equation}
where $\al^\pm_j, \be^\pm_j$ are both eigenforms of $\Delta_B$ 
with eigenvalue $\lambda_j^2$. We have proved
\begin{proposition} The indicial set $\Lambda$ for the
operator $D$ with respect to a $b$-metric $g$ consists of
the values $\pm \lambda$, where $\lambda^2 \in \mbox{spec}\,(\Delta_B)$
acting on $\left.{}^b\Omega^*(M)\right|_B$.
\label{pr:b-indr}
\end{proposition}
Note here that these calculations seem to leave open the possibility
that $0$ is a double root, which would allow for the possibility
of solutions of the indicial equation of the form $\omega = \omega' \log x
+ \omega'' x^0$. However, (\ref{eq:b9.5}) admits no solutions of 
this form, and so we see that the double root is spurious and 
arises merely from the algebraic calculations above. 

We conclude this section by discussing some general mapping properties
of $b$ operators on weighted $L^2$ spaces as well as regularity results
for their solutions.  Proofs of these theorems may be found in
\cite{Me-aps}.

Let $L^2_b(M) = L^2(M,\frac{dxdy}{x})$; this is the 
same as $L^2(M,dV_g)$ if $g$ is any $b$ metric. We also define
\[
H^\ell_b(M) = \{u \in L^2_b(M): V_1 \cdots V_j u \in L^2_b(M)
\ \forall \, j \leq \ell\ \mbox{and}\ V_i \in \calV_b\},
\]
and
\[
x^\gamma H^\ell_b(M) = \{u = x^\gamma v: v \in H^\ell_b(M)\},
\]
whenever $\ell \in \NN$ and $\gamma \in \RR$.

\begin{proposition} Let $P$ be an elliptic differential $b$ operator
of order $m$, acting between sections of the vector bundles $E$ and
$F$ over $M$, with indicial set $\Lambda$. Then the mapping
\[
P: x^\gamma H^{\ell+m}_b(M;E) \longrightarrow x^\gamma H^\ell_b(M;F)
\]
is Fredholm if and only if $\gamma \notin \{\mbox{Re}(\zeta):
\zeta \in \Lambda\}$.
\label{pr:mapb1}
\end{proposition} 

To state the final proposition, we introduce the important
subspace of polyhomogeneous distributions, sitting in the
space of conormal distributions:
\[
\calA^*_\phg(M) = \big\{u \in \calA^*(M): u \sim \sum_{{\mathrm{Re}}\, 
\gamma_j  \to \infty} \sum_{k = 0}^{N_j} u_{jk}(y)x^{\gamma_j}
(\log x)^k\ u_{jk} \in \calC^\infty(\del M)\big\}.
\]
These expansions are meant in the standard asymptotic sense as $x \to 0$ 
and hold along with all derivatives. The superscript $*$ here may be 
replaced by an index set $I$ containing all pairs 
$(\gamma_j,k)$ which are allowed to appear in this expansion.

\begin{proposition}
If $u \in x^\gamma L^2_b(M;E)$ and $Pu = 0$, then $u \in \calA^I_\phg(M;E)$,
where $I$ is an index set derived from the index set 
$\Lambda$ for $P$ truncated below the weight $\gamma$.  If $Pu =f$ where
$u \in x^\gamma L^2_b(M;E)$ and $f \in \calA^{\gamma'}(M;F)$ for
some $\gamma' > \gamma$, $\gamma' \notin \mbox{Re}\Lambda$,
then $u = v + w$ where $v \in \calA^I_\phg(M;E)$
and $w \in \calA^{\gamma'}(M;F)$.
\label{pr:mapb2}
\end{proposition}

The powers $\gamma$ appearing in the polyhomogeneous expansion in this
proposition are of the form $\gamma_j + \ell$ where each $\gamma_j$ is 
an element of the index set for $P$ and $\ell \in {\mathbb N}_0$. 
Logarithms can arise either from indicial roots with multiplicity greater 
than $1$, or else (as in classical ODE theory) when two indicial roots differ 
by an integer. For more details on this, see \cite{Me-aps}. All the roots 
we encounter in this paper are of multiplicity one 1 (although this
fact does not really affect the arguments much), and we shall justify
this in the various cases below, as we did following Proposition~\ref{pr:b-indr}.
 
These results about $b$ operators may be proved in a variety of
ways, some fairly elementary. For example, see \cite{APS} for 
the analysis of $\Delta_g$ on cylinders using separation of variables. 
We refer, however, to \cite{Me-aps} and \cite{ma-edge} for proofs based on
the calculus of $b$ pseudodifferential operators. This general theory is 
quite flexible, and is ideally suited for the proofs of more general 
index theorems in the $b$ category. A thorough treatment of this calculus, 
along with many applications, is given in \cite{Me-aps}.  

We shall not need to know much about these operators beyond
their mapping properties, but for
the sake of completeness, we say a few words about them. The calculus
${}^b\Psi^*(M)$ is designed in part to contain parametrices for
elliptic $b$-operators. Elements $A \in \Psi_b^*(M)$ are characterized
in terms of the structure of their Schwartz kernels $\kappa_A$.
Each such $\kappa_A$ is a distribution on $M^2 = M \times M$ with
singularities along the diagonal and side faces of this double space;
kernels of elements in ${}^b\Psi^*(M)$ are characterized by the fact
that they lift to distributions on a resolution $M^2_b$ of $M^2$
with only polyhomogeneous singularities. This resolution is the normal
blowup of $M^2$ along its corner and is obtained by replacing the corner
$(\del M)^2$ by its interior normal spherical bundle.

\subsection{Analysis for scattering metrics and operators}
We next consider scattering metrics on $M$. The analysis of general
elliptic operators in the scattering calculus is considerably
more subtle than for operators in the $b$ calculus, but because
we only consider the Hodge-de Rham operator, various simplifications
permit us to reduce directly to the $b$ calculus.
(Later in the paper, however, we shall need to use the calculus
of fibred boundary pseudodifferential operators, which is much
closer in spirit to the scattering calculus than to the
$b$ calculus.) 

Recall that a scattering metric $g$ has the form $g = g'/x^2$, where $g'$
is a $b$ metric. We define the Lie algebra $\calV_\sc$ of scattering
vector fields to consist of all smooth vector fields on
$\olM$ which have bounded length with respect to any scattering metric
$g$. Clearly
\[
\calV_\sc = x\calV_b = \{V: V = xW, \ W \in \calV_b\};
\]
alternately, in local coordinates $(x,y_1,\cdots, y_{n-1})$ near $\del M$,
$\calV_\sc$ is spanned by the vector fields $x^2\del_x$ and $x\del_{y_j}$.
By definition, these form the full set of sections of the scattering tangent
bundle ${}^\sc TM$; its dual, ${}^\sc T^*M$,
is locally smoothly trivialized by the sections
\[
\frac{dx}{x^2}, \frac{dy_1}{x}, \ldots, \frac{dy_{n-1}}{x}.
\]
The space of smooth sections of the exterior powers of this bundle is
$\calC^\infty\Omega^*_\sc(M)$. Thus any $\omega \in \calC^\infty
\Omega^*_\sc(M)$ can be written as
\begin{equation}
\omega = \sum_k \omega_k =
\sum_k \left(\frac{\al_k}{x^k} + \frac{dx}{x^2}\wedge
\frac{\be_{k-1}}{x^{k-1}}\right), \qquad
\al_k, \be_{k-1} \in \calC^\infty. 
\label{eq:scform}
\end{equation}
An advantage of this normalization is that
\[
||\omega||^2 = \int_M \sum_k \left(|\al_k|^2 + |\be_{k-1}|^2
\right) \, \frac{dxdy}{x^{n+1}}.
\]

An operator $P$ is a scattering differential operator if it can be
locally written as a finite sum of multiples of elements of $\calV_\sc$:
\[
P = \sum_{j+|\al|\leq m} a_{j,\al}(x,y)(x^2\del_x)^j(x\del_y)^\al,
\qquad a_{j,\al} \in \calC^\infty(\olM). 
\]
Its scattering symbol is defined as
\[
{}^\sc\sigma_m(P)(x,y;\xi,\eta) = i^{-m}
\sum_{j+|\al|=m} a_{j,\al}(x,y)\xi^j \eta^\al;
\]
$P$ is elliptic in this calculus if this symbol is invertible
for $(\xi,\eta) \neq 0$.

The analysis of $(\Delta - \lambda)u = 0$ is quite different
when $\lambda$ is negative or positive; for example, in the former case,
solutions decay rapidly while in the latter they oscillate with slow decay
as $x \to 0$. Accordingly, the nature of the resolvent changes
dramatically when $\lambda \in \mbox{spec}\,(\Delta)$ cf.\ \cite{Me-scm},
\cite{HV}. Because of this, the general theory of parametrices, mapping 
properties and regularity theory for elliptic scattering operators is 
fairly complicated. Fortunately we can sidestep this calculus by virtue 
of the 
\begin{proposition}
If $g$ is a scattering metric on $M$, then
\[
D = d + \delta: \calC^\infty\Omega_\sc^*(M) \longrightarrow
x \calC^\infty\Omega_\sc^*(M)
\]
is an elliptic first order scattering operator of the form $D = x D'$ where
$D'$ is an elliptic first order $b$-operator.
\end{proposition}
\begin{remark}
It seems initially somewhat confusing that $D'$ is a $b$-operator when
acting between sections of the scattering form bundles (normalized as
above), but not when acting between sections of the $b$ form bundles.  We
can understand why this is true, however, when we consider that the
endomorphism $dx \wedge$ has the same operator norm on forms as the does
the form $dx$. This norm depends upon the metric on
$M$.  Thus $dx/x$ is a unit norm endomorphism on the
bundle of forms when $M$ has a $b$-metric, whereas 
$dx/x^2$ is the unit endomorphism on the bundle of forms when $M$ has a
scattering metric.  There is a similar shift in the power of $x$ in the
other coordinates, so in the scattering case, an extra power of
$x$ is absorbed into the denominator of the endomorphism part of the
Laplacian.  This makes
$D'$ a $b$-operator on the bundle of scattering forms although it is not
as an operator on the bundle of $b$-forms.  
\end{remark}

\begin{proof} Write $\omega \in \calC^\infty\Omega^*_\sc(M)$
as in (\ref{eq:scform})
and set $\al = \sum \al_k$, $\be = \sum \be_k$. Then a brief
calculation gives
\begin{equation}
D\omega = \sum_k \left(
\frac{x(D_B\al)_k- x^2 \del_x \be_k +
(n-k-1)x \be_k}{x^k} +
\frac{dx}{x^2}\wedge
\frac{x^2\del_x \al_k - kx \al_k - x(D_B \be)_{k}}
{x^k}\right),
\label{eq:dsc2}
\end{equation}
where $(D_B \zeta)_k$ is the component of degree $k$ of $D_B \zeta$ for
$\zeta = \al$ or $\be$. This shows immediately that $D' \equiv x^{-1}D$
is a $b$ operator; it differs from $D_{g'}$, where $g' = x^2 g$ is the
associated $b$ metric, only in terms of order zero.  Of course, these
affect the indicial set $\Lambda$ markedly. 
\end{proof}

The mapping properties of $D$ and the regularity properties of its
solutions may be deduced directly from the corresponding properties
for $D'$ in Propositions \ref{pr:mapb1} and \ref{pr:mapb2}.
Note, however, that the extra factor of $x$ causes a shift in the
weight of the function spaces.

\begin{proposition}
Suppose $g$ is a scattering metric, so that $D = xD'$ as above.
Let $\Lambda$ denote the indicial set for $D'$. Then
\[
D: x^{\gamma}H^{\ell+1}_b\Omega^*_{\sc}(M) \longrightarrow
x^{\gamma+1}H^{\ell}_b\Omega^*_{\sc}(M)
\]
is Fredholm for any $\ell \in \NN_0$ and $\gamma \notin \{\Re(\lambda):
\lambda \in \Lambda\}$. 
\label{pr:mapsc1}
\end{proposition}
\begin{proposition} If $\omega \in x^\gamma L^2\Omega^*_\sc(M)$ 
for any $\gamma \in \RR$ and
$D\omega = 0$, then $\omega \in \calA^I_\phg\Omega^*_\sc(M)$, where
$I$ is some augmented index set determined by the indicial set $\Lambda$
of $D'$ and the cutoff weight $\gamma$. If, on the other hand,
$D\omega = \eta$ where $\eta \in x^{\gamma'+1}\calA^*\Omega^*_\sc(M)$ for
some
$\gamma' > \gamma$, then $\omega = \omega' + \omega''$ with $\omega'
\in \calA^*_\phg\Omega^*_\sc(M)$ and $\omega'' \in
\calA^{a+1}\Omega^*_\sc(M)$.
\label{pr:mapsc2}
\end{proposition}

We conclude this section with a computation of the relevant part 
of the indicial set $\Lambda$ for $D'$. As in the $b$ case, this 
set is determined by the spectrum of $\Delta_B$, but the computation is 
more intricate. 

First, define the numerical operators $N_1$ and $N_2$
\[
N_1 \be_k = (n-k-1)\be_k, \qquad N_2 \al_k = -k \al_k
\]
(i.e.\ $N_1$ and $N_2$ are diagonal on $\Omega^*(B)$ with respect
to the decomposition by degree). Let
\[
\omega = \sum \omega_k, \qquad
\omega_k = \frac{\al_k}{x^k} + \frac{dx}{x^2}\wedge
\frac{\be_{k-1}}{x^{k-1}},
\]
where all $\al_j$ and $\be_j$ are independent of $dx$. Then
$D(x^\gamma \omega) = x^{\gamma+1}I_{D'}(\gamma)(\omega)$ where
\[
I_{D'}(\gamma)(\omega) =
\sum_k \left(\frac{(D_B \al + (N_1 - \gamma)\be)_k}{x^k}
+ \frac{dx}{x^2}\wedge \frac{(-D_B\be +
(N_2 + \gamma)\al)_{k-1}}{x^{k-1}}\right).
\]
Writing $I$ for $I_{D'}$, this vanishes when
\begin{equation}
I(\gamma)
\left(\begin{array}{c}
\al \\ \be
\end{array}
\right) \equiv 
\left(
\begin{array}{cc} 
D_B & N_1 - \gamma \\ N_2 + \gamma & -D_B
\end{array} 
\right)
\left(\begin{array}{c}
\al \\ \be
\end{array}
\right) = 
\left(
\begin{array}{c}
0 \\ 0
\end{array}
\right)
\label{eq:scind-r}
\end{equation}

Although this equation seems strongly coupled, and hence difficult
to analyze, computations can be simplified using the special structure
that $\Delta = D^2$ preserves degree. On the indicial level this gives
\[
I(\gamma+1)I(\gamma)
\left(\begin{array}{c}
\al \\ \be
\end{array}
\right) =
\left(
\begin{array}{cc} 
D_B & N_1 - \gamma -1 \\ N_2 + \gamma +1 & -D_B
\end{array} 
\right)
\left(
\begin{array}{cc} 
D_B & N_1 - \gamma \\ N_2 + \gamma & -D_B
\end{array} 
\right)
\left(\begin{array}{c}
\al \\ \be
\end{array}
\right) = 
\left(
\begin{array}{c}
0 \\ 0
\end{array}
\right).
\]
Multiplying out this matrix of operators and using the easily
verified fact that
\[
[D_B,N_j] = d_B - \delta_B,
\]
we have
\[
\left(
\begin{array}{cc} 
\Delta_B + (N_1-\gamma-1)(N_2 + \gamma) & 2d_B \\
2\delta_B & \Delta_B + (N_2+\gamma+1)(N_1 - \gamma)
\end{array} 
\right)
\left(
\begin{array}{c}
\alpha \\ \beta
\end{array}
\right)
= \left(
\begin{array}{c}
0 \\ 0
\end{array}
\right).
\]
The coupling here occurs only between closed $k$ forms
and coclosed $(k-1)$-forms.

We shall not need to calculate all the indicial roots of $D'$,
although this can be done readily from these formul\ae.
Instead, we focus on the special value $\gamma = n/2 - 1$.
This is a critical value in our calculations because $x^{n/2}$
is just on the border of lying in $L^2(dV_g) = L^2(x^{-n-1}dxdy)$
and we shall need to analyze the map $D: x^{\gamma-1}L^2 \to
x^\gamma L^2$ for $\gamma$ near this borderline value.
Thus, setting $\gamma = n/2-1$ gives
\begin{eqnarray*}
(N_1- n/2)(N_2 + n/2-1)\al_k = (n/2 - k -1)^2 \al_k, \\
(N_2 + n/2)(N_2 - n/2+1)\be_{k-1} = (n/2 - k + 1)^2 \be_{k-1}.
\end{eqnarray*}
Hence, if $\omega$ lies in the nullspace of $I_D(n/2 - 1)$ then
for all $k$ we have:
\begin{eqnarray*}
(\Delta_B  + (n/2 - 1 - k)^2) \al_k  + 2d_B\be_{k-1} = 0 \\
(\Delta_B  + (n/2 + 1 - k)^2) \be_{k-1}  + 2\delta_B\al_{k} = 0 .
\end{eqnarray*}
Decompose these equations using an eigendecomposition for $\Delta_B$ such
that $\al_k = a \psi_k,\ \be_{k-1} = b \phi_{k-1}$,
where both $\psi_k$ and $\phi_{k-1}$ are eigenforms with
eigenvalue $\lambda^2 \geq 0$ and $d\phi_{k-1} = \lambda \psi_k$,
$\delta \psi_k = \lambda \phi_{k-1}$. Then
\[
\left(
\begin{array}{cc} 
\lambda^2 + (n/2 - k - 1)^2 & 2\lambda \\
2\lambda & \lambda^2 + (n/2 - k + 1)^2
\end{array} 
\right)
\left(
\begin{array}{c}
a \\ b
\end{array}
\right)
= 
\left(
\begin{array}{c}
0 \\ 0
\end{array}
\right),
\] 
and so there are nontrivial solutions only if this matrix is singular.
Its determinant equals $(\lambda^2 + (n/2 - k)^2 - 1)^2$, hence there
are no solutions unless $|k-n/2|\leq 1$. First, if $\lambda = 0$, then
$k = n/2 \pm 1$, and the nullspace consists of harmonic forms
$\alpha_{n/2 - 1}$ and $\beta_{n/2 + 1}$.  Next, if
$\lambda^2 = 1$ is in the spectrum of $\Delta_B$, then there are solutions
for $\alpha_k$ and $\beta_{k-1}$ only if $k = n/2$, and elements
of the nullspace of $I(n/2)I(n/2-1)$ are obtained by taking $a +b = 0$. 
Finally, there are solutions of a similar type when $k = (n \pm 1)/2$
and $\lambda^2 = 3/4 \in \mbox{spec}\,(\Delta_B)$.  One can then
verify that only the solutions corresponding to $\lambda=0$
also lie in the nullspace of $I(n/2-1)$, hence these are the only ones
which appear in the polyhomogeneous expansions for solutions of
$D\omega = 0$. 

Note that $\gamma = n/2 - 1$ is not an indicial root of multiplicity
two. As in the $b$ setting, this follows by checking that there are
no solutions of (\ref{eq:scind-r}) of the form $\omega' \log x + 
\omega''$. 

We conclude these computations by noting that if $\omega \in L^2(dV_g)$
satisfies $D\omega = 0$, then $D^2 \omega = 0$ and the usual integration
by parts, which is justified in $L^2$, gives $d\omega = \delta \omega = 0$
individually. To relate this to the preceding calculations, this implies
that if $\gamma > n/2$ is an indicial root for $x^{-1}D$, then it is
also one for both $x^{-1}d$ and $x^{-1}\delta$ (and conversely),
and these are much simpler to compute. In fact,
\[
I_{x^{-1}d}(\gamma)(\omega) = \sum_k \left(
\frac{d_B \al_k}{x^{k+1}} + \frac{dx}{x} \wedge
\frac{(\gamma-k)\alpha_k - d_B \be_{k-1})}{x^k} \right)
\]
and
\[
I_{x^{-1}\delta}(\gamma)(\omega) = \sum_k \left(
\frac{\delta_B \al_k + (n-k-\gamma)\be_{k-1}}{x^{k-1}} +
\frac{dx}{x}\wedge \frac{-\delta_B\be_{k-1}}{x^{k-2}}\right).
\]
Hence $\omega$ is in the nullspace of both these operators
provided $d_B \al = \de_B \be = 0$ and also
\[
\delta_B \al_k = -(n-k-\gamma)\be_{k-1}, \
d_B \be_{k-1} =(\gamma-k)\al_k.
\]
Thus $\al_k$ is closed, $\be_{k-1}$ is coclosed, and both are in the
nullspace of $\Delta_B + (\gamma-k)(n-k-\gamma)$. On an eigenspace with
eigenvalue $\lambda^2$ for $\Delta_B$ we must have
\[
\gamma^2 - n\gamma + (k(n-k) - \lambda^2) = 0,
\]
and by assumption above we must choose the root which
is greater than $n/2$. (Of course, solutions of these
equations, no matter the value of $\gamma$, give indicial
roots of $x^{-1}D$, corresponding to non-$L^2$ solutions. The
point of the earlier calculations is that there are other
indicial roots, corresponding to solutions which are
not individually closed and coclosed.) In summary, these comprise
the subset 
\begin{equation}
\Lambda' = \{\gamma_j^\pm: \mbox{roots of}\ \gamma^2 - n\gamma +
k(n-k) - \lambda_j^2 =0,\
\lambda_j^2 \in \mbox{spec}(\Delta_B)\}
\label{eq:indsetsc}
\end{equation}
inside the possibly larger set of all indicial roots of $x^{-1}D$.
Note in particular that when $\lambda_j^2 = 0$, $\gamma_j^\pm =
k, n-k$.  

\subsection{Hodge theorems for $b$ and scattering metrics}

Having assembled these analytic facts and calculations, we now
complete the proofs of the Hodge theorems for $b$ and scattering
metrics following the outline from the compact case. We invert
the usual order of presentation and discuss first the $b$ case, which
specializes Theorem 2, and afterwards the scattering case,
which specializes Theorem 1. These results equate Hodge cohomology
with weighted cohomology only; Corollary \ref{cor:luxembourg} 
shows that the results are indeed the same as stated in
Theorems 2A and 1A, respectively. 

\medskip
\noindent{\bf Theorem 2B.}\ {\it
Let $g$ be an exact $b$ metric on the manifold $M$. Then for
sufficiently small $\e > 0$ and for every $k = 0, \ldots, n$, 
there is a canonical isomorphism
\begin{equation}
\Phi: L^2\calH^k(M) \longrightarrow \mbox{\rm Im}\, (WH^k(M,g,\e)
\longrightarrow WH^k(M,g,-\e)).
\label{eq:bPhi}
\end{equation}
}

\begin{proof}
As in the compact case, if $\omega \in L^2\calH^k(M)$, then $d\omega=0$.
Further, by Proposition \ref{pr:mapb2} and (\ref{eq:b13}),
$\omega$ is polyhomogeneous, with an
expansion of the form $\sum \omega_j^\pm(y)x^{\pm \lambda_j}$, where
the tangential and normal parts of $\omega_j^\pm$ are eigenforms on
$\del M$ with eigenvalue $\lambda_j^2$. Since $\omega \in L^2$, we see
that all coefficient forms $\omega_j^-$ vanish, as do those
$\omega_j^+ = 0$ corresponding to values of $j$ with $\lambda_j = 0$.
Hence $\omega = \alpha + \frac{dx}{x}\wedge \beta$ where
$\alpha,\beta = {\mathcal O}(x^{\underline{\lambda}})$,
where $\underline{\lambda} = \inf\{|\lambda_j| \neq 0: \lambda^2_j
\in \mbox{spec}(\Delta_B)\}$. Thus $[\omega] \in WH^k(M,g,\e)$ is
well-defined provided $\e < \underline{\lambda}$.

If $\omega \in L^2\calH^k(M)$ and $\Phi(\omega)=0$, then
$\omega = d\zeta$ for some $\zeta \in x^{-\e} L^2\Omega_b^{k-1}(M)$.
Computing cohomology with the complex of conormal forms as
explained in \S 2, we can take $\zeta = \mu + \frac{dx}{x}\wedge \nu$,
where $\mu, \nu \in \calA^{-\e}([0,1)_x\times B; \bigwedge^*(B))$.
In the integration by parts $||\omega||^2 = <\omega, d \zeta> =
<\delta \omega,\zeta>= 0$, the boundary term equals $\lim_{x \rightarrow 0}
<\alpha, \nu>_B$, and this vanishes since $\e < \underline{\lambda}$.
Hence $\omega = 0$ and so $\Phi$ is injective.

Next, by Proposition \ref{pr:mapb1}, the map
\begin{equation}
D: x^{-\e}H^1_b\Omega^*(M) \longrightarrow x^{-\e}L^2\Omega^*(M)
\label{eq:FDb}
\end{equation}
is Fredholm when $\e \in (0,\underline{\lambda})$. This gives the
decomposition
\[
x^{-\e}L^2\Omega^*(M,dV_g) = \left(\mbox{ran}\, \left. D
\right|_{x^{-\e}H^1_b\Omega^*}
\right) \oplus \left(\mbox{ran}\, \left. D
\right|_{x^{-\e}H^1_b\Omega^*}\right)^\perp.
\]
The second summand on the right is finite dimensional, and could
be replaced with any other finite dimensional subspace of $x^{-\e}
L^2\Omega^*$ which is complementary to the range of $D$, since the
orthogonality of this decomposition does not play any role. In 
particular, we claim that we can replace this term by $L^2\calH^*(M)$.
To see this, note simply that the natural pairing between $x^{-\e}L^2\Omega^*$ 
and $x^{\e}L^2\Omega^*$ identifies the orthogonal complement of the 
range of $D$ on $x^{-\e}L^2\Omega^*$ with the nullspace of $D$ on 
$x^{\e}L^2\Omega^*$, which equals $L^2\calH^*(M)$. In any case,
we have shown that for any $\eta \in x^{-\e}L^2\Omega^*$, there exist 
elements $\zeta \in x^{-\e}H^1_b\Omega^*$ and $\gamma \in L^2\calH^*$ such that
\begin{equation}
\eta = D\zeta + \gamma.
\label{eq:eta}
\end{equation}
 

Now we prove surjectivity of $\Phi$. Fix any $[\eta]$ in the
space on the right in (\ref{eq:bPhi}) and choose a conormal
representative $\eta \in \calA^{\e}\Omega^*$ for
it. Decompose $\eta$ as $D\zeta + \gamma$ as above, in the space
$x^{-\e}L^2\Omega^*$. Proposition \ref{pr:mapb2} shows that
$\zeta$ is partially polyhomogeneous, i.e.\ it is a sum of a
finite number terms of
the form $\zeta_{j,\ell}\, x^{\sigma_j}(\log x)^{\ell}$ and a term
$\zeta' \in \calA^\e\Omega^*$. All exponents $\sigma_j$ lie in the
interval $(-\e,\e)$, and because we can choose $\e$ as small
as desired, we may assume that the only terms which appear have
$\sigma_j = 0$. The remaining terms correspond to solutions of the
indicial operator $I_D(0)$, and the analysis in \S 4.1 shows
that $0$ is an indicial root of multiplicity one, so no log terms
occur. Thus $\zeta = \zeta_0 + \zeta'$ where $I_D(0)\zeta_0 = 0$;
writing $\zeta_0 = \mu_0 + (dx/x)\wedge \nu_0$ then both $\mu_0$ and $\nu_0$
are harmonic on $B$.

As in \S 3, to conclude that $\delta \zeta = 0$
we must check that all three terms, $\langle \de \zeta,
\eta \rangle$, $\langle \de \zeta, \gamma \rangle$ and
$\langle \de \zeta, d\zeta \rangle$,  vanish. This is true formally,
i.e.\ integrating by parts and neglecting the boundary terms, so
it suffices to check that these boundary terms also vanish.
For the first two this is straightforward since $|\zeta|$ is bounded
and both $\eta$ and $\gamma$ vanish at $x=0$. For the final term, the
boundary contribution is
$$
\int_M d(\zeta \wedge * \delta \zeta)
= \int_B \mu_0 \wedge d_B *_B \nu_0 = <\mu_0, \delta_B \nu_0>_B,
$$
and this vanishes since $\nu_0 \in L^2\calH^*(B)$.
Hence $\de \zeta = 0$, and so $\eta = d\zeta +
\gamma$, as required. As in the compact case, there are only
forms of degree $k$ here. This finishes the proof. 
\end{proof}

Now suppose that
\[
g = \frac{dx^2}{x^4} + \frac{h}{x^2}
\]
is a scattering metric on $M$.

\medskip
\noindent{\bf Theorem 1B.}\ {\it
Let $g$ be a scattering metric on $M$. Then for any $\e>0$ sufficiently
small, there is a canonical isomorphism
\begin{equation}
\Phi:  L^2\calH^k(M) \longrightarrow
\mbox{\rm Im}\,(\calW H^k(M,g,\e) \to \calW H^k(M,g,-\e)).
\label{eq:scPhi}
\end{equation}
\label{th:htsc}  
}

\medskip

\begin{proof} If $\omega \in L^2\calH^k(M,g)$ then
$d\omega = 0$. By Proposition (\ref{pr:mapsc2}), $\omega
\in \calA^{\underline{\lambda}}\Omega^k_\sc(M)$ for
some $\underline{\lambda} > 0$. Hence $[\omega] \in \calW H^k(M,g,\e)$
is well-defined provided $0 < \e < \underline{\lambda}$. Thus
$\Phi(\omega)$ is well-defined.

Suppose $\Phi(\omega) = 0$, so that $\omega = d\zeta$ for some
$\zeta \in x^{-\e-1}L^2\Omega_\sc^{k-1}$. Write
\[
\omega = \frac{\al}{x^k} + \frac{dx}{x^2}\wedge \frac{\be}{x^{k-1}},\qquad
\mbox{and} \quad \zeta = \frac{\mu}{x^{k-1}} + \frac{dx}{x^2}\wedge
\frac{\nu}{x^{k-2}}.
\]
We may assume that $\zeta$ is conormal, and hence that $|\mu|, \, |\nu|
\in {\mathcal O}(x^{n/2 + \e' - 1})$ for some $\e'>\e$. This implies
that $\lim_{x \rightarrow 0} <x^{-k+1} \mu , x^{-n+k} \beta>_B = 0$,
whence $||\omega||^2 = <d\zeta,\omega> = < \zeta,\delta \omega > = 0$.
This shows that $\omega = 0$ and thus $\Phi$ is injective.

The surjectivity argument proceeds as before. Since
\[
D: x^{-\e-1}H^1_b\Omega^*(M) \longrightarrow x^{-\e}L^2\Omega^*(M)
\]
is Fredholm, we have
\[
x^{-\e}L^2\Omega^*(M) = \left( \left. \mbox{ran}\, D
\right|_{x^{-\e-1}H^1_b}
\right) \oplus 
\left( \left. \mbox{ran}\, D \right|_{x^{-\e-1}H^1_b}
\right)^\perp.
\]
The same argument as in the $b$ case identifies this orthocomplement with
$L^2\calH^*$.  Now let $\eta \in \calA^\e\Omega^*$ represent a
nontrivial class $[\eta]$ in the space on the right in (\ref{eq:scPhi}).
Write $\eta = D\zeta + \gamma$, $\gamma \in L^2\calH^*$; by
Proposition \ref{pr:mapsc2}, $\zeta$ is partially polyhomogeneous
and is a finite sum of terms $\zeta_{j,\ell}x^{\sigma_j}(\log
x)^\ell$ and some $\zeta' \in \calA^\e\Omega_\sc^*$. By taking
$\e$ small enough, we can eliminate all but the term of weight
$n/2-1$, and by the computations in \S 4.3, this indicial root
occurs with multiplicity one so there are no log terms. In fact, those
computations give that
\[
\zeta = \zeta_0 + \zeta',\quad \zeta_0 = \left(\frac{\al_{n/2-1}}{x^{n/2 -1}} +
\frac{dx}{x^2}\wedge \frac{\be_{n/2+1}}{x^{n/2+1}}\right)x^{n/2 -1}
+ \zeta', \quad |\zeta'| = {\mathcal O}(x^{\e'}),
\]
where $\alpha_{n/2-1}$ and $\beta_{n/2+1}$ are harmonic on $B$.
Using the same reasoning as in the previous proof, we see that the
boundary terms in the integrations by parts $\langle \delta \zeta,
\eta\rangle = \langle \zeta, d\eta \rangle$
and $\langle \delta \zeta, \gamma \rangle = \langle \zeta, d\gamma \rangle$
both vanish. Since $d\eta = d\gamma = 0$, these terms vanish altogether.
Finally, $\langle \delta \zeta, d\zeta \rangle$ equals the sum
of the (vanishing) interior term, $\langle \zeta, d^2\zeta \rangle = 0$,
and a boundary term. Since $d$ and $\delta$ are both $x$ times
$b$-operators, 
$d\zeta'$ and $\delta \zeta'$ both decay, so only $\zeta_0$ contributes.
This boundary term equals
\[
\int_M d\zeta_0 \wedge * \delta \zeta_0 = \pm
\int_M d(\zeta_0 \wedge d * \zeta_0) = \pm
\langle \alpha_{n/2-1}, \delta_B\beta_{n/2+1}\rangle = 0.
\]
Hence $\delta \zeta = 0$, and finally, $\eta = d\zeta + \gamma$
where $\gamma \in x^{-\e-1}L^2\Omega^{k-1}$ and $\gamma \in L^2\calH^k$.
\end{proof}

\section{Fibred ends}
We now turn to the general case where both the base and fiber in
the fibration of $\del M$ are nontrivial. As in the $b$- and
scattering cases, we must: determine the explicit structure of $D$,
calculate its indicial roots, and understand its mapping properties
and the regularity (polyhomogeneity) of elements of $L^2\calH^*$.
For the construction of a parametrix for $D$, we invoke
the fibred boundary calculus of pseudodifferential operators,
as developed in \cite{MaMe} and extended in \cite{Va}. This serves as
a replacement for the $b$-calculus in this context, but
is more intricate. To help mitigate the analytic requisites, we
include a discussion of this parametrix construction in the very special
case where $M$ is a global product and the fibred boundary or fibred cusp
metric respects this  decomposition. Although the Hodge theorems in
these cases follow directly via a K\"unneth theorem from those for
$b$- and scattering metrics, we sketch an explicit parametrix
construction for $D$ in hopes that this gives some insight into
the more general case.

We begin with a general discussion of the fibred boundary
calculus, and then proceed immediately to a discussion of
$D$ and its parametrix in the product case. This is followed
by a review of the geometry of fibrations
and the structure of $D$ in the general case. The identifications
of Hodge cohomology with weighted cohomology are then proved,
as usual following the general line of argument from \S 3.
The final subsection relates the weighted cohomology to
intersection cohomology.

\subsection{The fibred boundary calculus}
Suppose that $\phi: Y=\del M \to B$ is a fibration with fiber $F$,
$\dim B = b$ and $\dim F = f$. Fixing an extension of this
fibration to a collar neighbourhood $\calU$ of $\del M$ in $M$,
we choose a fibred boundary metric
\[
g_\fb = \frac{dx^2}{x^4} + \frac{\phi^*(h)}{x^2} + k_F,
\]
where $h$ is a metric lifted from $B$ and $k_F$ is a symmetric
$2$-tensor which restricts to a metric on each fibre. There is an
associated fibred cusp metric $g_\fc = x^2 g_\fb$. 
These metrics stand in the same relationship to one another as do
scattering and $b$ metrics. As we saw in those cases, only the
$b$ calculus (but not the scattering calculus) is required to
analyze the operator $D$ in both cases. Similarly, the fibred
boundary calculus is enough to analyze the Hodge-de Rham operators
for both $g_\fb$ and $g_\fc$. (Indeed, there is no calculus
directly associated to $g_\fc$, for reasons indicated below.)

The fibred boundary calculus relies on the choice of a $1$-jet
of the defining function $x$ along the fibres at $\del M$. Making 
such a choice, define the Lie algebra of fibred boundary vector fields 
\[
\calV_\fb = \{V \in \calV_b(M): V \ \mbox{tangent to fibres
$F$ at $\del M$},\ Vx = {\mathcal O}(x^2)\}.
\]
To understand this more clearly, choose local coordinates $(x,y,z)$ 
where $y$ are coordinates on $B$, pulled back to $Y$ via $\phi$ and 
then extended into the manifold, $z$ are functions on $Y$ which restrict 
to coordinates on the fibres, similarly extended inward, and 
$x$ is in the given equivalence class of defining functions. Then $\calV_\fb$
is spanned locally over $\calC^\infty$ by the vector fields
$x^2\del_x, x\del_{y_j}, \del_{z_\ell}$. If $(\tilde{x},\tilde{y},\tilde{z})$ 
is a new choice of coordinates adapted to the fibration, then
$\del_{\tilde{z}}$ transforms into a vector field with one component equal to
$(\del x/\del \tilde{z}) \del_{x}$, and this explains why we need
to fix the differential of $x$ along each fiber, in order that the
coefficient here vanish to second order. 

In contrast, the vector fields 
associated to a fibred cusp metric are $x\del_x$, $\del_{y_j}$,
$x^{-1}\del_{z_\ell}$; these are singular,
but much more seriously, their span is not closed under Lie bracket.
Involutivity is a basic requirement in the microlocalization
procedure leading to the construction of the associated pseudodifferential
calculus, and this explains why there is no separate fibred cusp calculus.
The elements of $\calV_\fb$ constitute the full set of sections of the
$\fb$ tangent bundle ${}^\fb TM$. We use its dual, the $\fb$
cotangent bundle, and the bundle of $\fb$ exterior forms,
$\bigwedge^*_\fb(M)$. Thus, in the coordinates above,
\[
\calC^\infty\Omega^k_\fb(M) \ni \omega = \sum_{i=0}^k \frac{\al_{i}}{x^i} +
\frac{dx}{x^2}\wedge \sum_{j=0}^{k-1} \frac{\be_j}{x^j},
\]
where $\al_i$ is a sum of wedge products of $i$-forms in $y$ and
$k-i$ forms in $z$ and $\be_j$ is a sum of wedge products of
$j$-forms in $y$ and $k-j-1$ forms in $z$, all of which are smooth
in the ordinary sense on $\olM$. (This decomposition is recast
more invariantly later.)

We now define the space of $\fb$ differential operators on $M$,
the associated $\phi$ symbol, and finally, the corresponding notion
of symbol ellipticity. This leads to the
\begin{proposition} For an exact fibred boundary metric, the Hodge-deRham
operator $D = d+\delta$ is an elliptic first order fibred boundary
differential operator. For an exact fibred cusp metric, the operator
$D$ is of the form $x^{-1}D'$ where $D'$ is an elliptic first order
fibred boundary operator.
\end{proposition}

Elliptic fibred boundary operators may be analyzed using the
calculus of fibred boundary pseudodifferential operators from
\cite{MaMe}. We shall use the elaboration of this theory developed
by Vaillant \cite{Va}. He constructs parametrices for any Dirac-type
operator associated to a fibred boundary or fibred cusp metric, and
in particular proves
\begin{proposition}[\cite{Va}, Proposition 3.28]
Let ${\mathcal D}$ be a Dirac-type operator associated to a fibred
boundary metric (for example, either $D$ or $D'$ above). Suppose that
$\omega \in x^\gamma L^2\Omega^*_\fb(M)$ satisfies
${\mathcal D}\omega = 0$. Then $\omega \in \calA^*_\phg\Omega^*_\fb(M)$.
\label{pr:Vaphg}
\end{proposition}

We also require a replacement for the other parts of Propositions
\ref{pr:mapb2} and \ref{pr:mapsc2}, as well as replacements for
the basic mapping properties, as in Propositions \ref{pr:mapb1}
and \ref{pr:mapsc1}. The precise forms of these results in the
fibred boundary setting are somewhat different, and as explained
in the preamble to this section, to motivate these results we shall 
take a detour and investigate the mapping and regularity 
properties for product metrics. This involves little more
than rephrasing the corresponding results for $b$ and scattering
metrics, but is included to help orient the reader. We also
include a discussion of the indicial root structure for $D$
in these two cases; the computations are more transparent
in the product cases, but the general results are qualitatively the same.

\subsection{The product case}
\subsubsection{Fibred boundary metrics}
Suppose that $M = N \times F$, where $\del N = B$, and
fix a fibred boundary metric $g$ on $M$ which is of the form
$g' + k$ where $g'$ is a scattering metric on $N$ and $k$ is a
metric on the compact manifold $F$.

We write $Y = \del M = B \times F$. Since $TY$ splits canonically
as $TB \oplus TF$, we have
\[
\mbox{$\bigwedge$}^kT^*Y = \bigoplus_{p+q=k}
\mbox{$\bigwedge$}^{p,q} Y, \qquad \mbox{$\bigwedge$}^{p,q}Y =
\mbox{$\bigwedge$}^p T^*B \otimes \mbox{$\bigwedge$}^q T^*F.
\]
Thus any $\omega \in \bigwedge^k M$ can be written as
\[
\omega = \frac{\alpha}{x^k} + \frac{dx}{x^2} \wedge \frac{\beta}{x^{k-1}},
\qquad \mbox{where}\qquad
\alpha \in \bigoplus_{j}\mbox{$\bigwedge$}^{k,j}Y, \quad \beta \in
\bigoplus_j \mbox{$\bigwedge$}^{k-1,j}Y
\]
depend parametrically on $x$.

The Hodge-de Rham operator $D = D_M$ acts on $\omega$, regarded as
a column vector $(\alpha,\beta)^t$, as
\begin{equation}
\left(
\begin{array}{cc} 
0 & -x^2 \del_x + (b-k+1) x \\
x^2 \del_x - k x & 0
\end{array} 
\right)
+ 
\left(\begin{array}{cc}
x D_B + D_F & 0 \\ 0 & -x D_B - D_F \end{array} \right).
\label{eq:DYprodfb}
\end{equation}
Here $D_F$ acts on a $(p,q)$ form $\eta \wedge \nu$ as $(-1)^p \eta
\wedge (D_F \nu)$. In the more general (nonproduct) case,
$D$ has a similar decomposition, but the second matrix
has extra terms coming from the nontrivial geometry of the bundle.

The space of harmonic forms on the compact manifold $F$ is finite
dimensional. Let 
\[
\Pi_0: L^2 \Omega^* (F) \longrightarrow L^2 \calH^*(F),
\qquad \Pi_\perp = I - \Pi_0
\]
be the natural orthogonal projectors. These extend naturally
to $L^2\Omega_\fb^*(M)$, and we have
\[
D_M = \Pi_0 D_M \Pi_0 + \Pi_\perp D_M \Pi_0 + \Pi_0 D_M \Pi_\perp
+ \Pi_\perp D_M \Pi_\perp.
\]
Since $[D_M,\Pi_0]=0$ in the product case, the second and third terms
vanish and this reduces to
\[
D_M=\Pi_0 D_M \Pi_0 \oplus \Pi_\perp D_M \Pi_{\perp}.
\]

We use this decomposition to construct a parametrix for $D_M$. First
\[ 
\Pi_0 D_M \Pi_0 = D_N \otimes \mbox{Id}_{\calH^*(F)},
\]
and so by the theory from \S 4.1 and 4.2, if $a \in \RR$
is not an indicial root for $x^{-1}D_N$, this operator
is Fredholm as a mapping from $x^{a}L^2\Omega^* \to x^{a+1}L^2\Omega^*$.
We write the generalized inverse as
\[
G_0^a: x^{a+1}L^2\Omega^*_\fb(N)\otimes\calH^*(F)
\longrightarrow x^a H^1\Omega^*_\fb(N)\otimes\calH^*(F).
\]

The second term in the decomposition of $D_M$ has square $\Delta_N + 
\Pi_\perp\Delta_F\Pi_\perp$. Since $\Delta_N \geq 0$ and $\Pi_\perp 
\Delta_F \Pi_\perp \geq c > 0$, we have that $\Pi_\perp D_M \Pi_\perp: 
x^a H^1\Omega^*_\fb \to x^a L^2\Omega^*_\fb$ is an isomorphism 
for any $a$. Thus for any $a$ we get 
\[
G^a_\perp \equiv (\Pi_\perp D_M \Pi_\perp)^{-1}: x^a 
\Pi_\perp L^2\Omega^*_\fb(M) \longrightarrow
x^a \Pi_\perp H^1 \Omega^*_\fb(M).
\]
Altogether, we have proved that
\[
G^a_0 \oplus G^a_\perp \equiv G^a: x^{a+1}\Pi_0 L^2\Omega^*_\fb(M) 
\oplus x^a\Pi_\perp L^2\Omega^*_\fb(M)
\longrightarrow x^a H^1\Omega^*_\fb(M)
\]
is bounded. Clearly $I - G^a D_M = \Pi_M^a$ is the projector onto
the nullspace of $D_M$ in $x^a L^2\Omega^*_\fb$, which is the
same as the nullspace of $\Pi_0 D_M \Pi_0$, i.e.\ $L^2\calH^*(N)
\otimes \calH^*(F)$. By Proposition \ref{pr:mapsc2},
$\Pi_M^a$ maps $x^a L^2$ into the space of polyhomogeneous fiber harmonic
forms. We emphasize that the indicial roots for $D_M$ are exactly
the same as for $D_N$. In particular, the `critical' root
$\frac12\dim N - 1 = (b-1)/2$ has multiplicity one!

In summary, we have proved that the mappings
\begin{equation}
D_M: x^a H^1\Omega^*_\fb(M) \longrightarrow x^{a+1}\Pi_0 L^2\Omega^*_\fb(M)
\oplus x^{a}\Pi_\perp L^2\Omega^*_\fb(M)
\label{eq:prm1}
\end{equation}
and
\begin{equation}
D_M: x^{a-1}\Pi_0 H^1\Omega^*_\fb(M) \oplus x^a \Pi_\perp H^1\Omega^*_\fb(M)
\longrightarrow x^a L^2\Omega^*_\fb(M)
\label{eq:prm2}
\end{equation}
are Fredholm when $a$, respectively $a-1$, is not an indicial
root of $D_N$. 

The generalized inverse $G^a$ has other mapping properties.
Suppose $\eta = D_M \zeta$ where $\eta \in \calA^a \Omega^*_\fb(M)$ and 
$\zeta \in x^{c-1}\Pi_0 H^1\Omega^*_\fb(M) \oplus x^{c} \Pi_\perp H^1_\fb(M)$ 
for some $c<a$. Then in fact $\zeta \in \Pi_0 \calA^*_{\phg}
\Omega^*_{\fb}(M) + \calA^{a}\Omega^*_{\fb}(M)$. 

\subsubsection{Fibred cusp metrics}
Now suppose that $M=N \times F$ has a fibred cusp metric
$g_\fc$; notice that this is a warped product (since
$k_F$ is multiplied by $x^2$). We obtain a parametrix for the
associated Hodge-de Rham operator $D$ as above.
Write all forms in terms of the (essentially
orthonormal) coframe, $dx/x$, $dy$, $x\, dz$, and denote
the space of forms with this normalization as $\bigwedge^*_\fc$. Thus
\[
\Lambda^*_\fc(M) \ni \omega = x^k \alpha + \frac{dx}{x}
\wedge x^k \beta, \qquad \mbox{where}\quad
\alpha, \beta \in \bigoplus_j \mbox{$\bigwedge$}^{j,k}Y.
\]
Write $\omega \in \Omega^{*,k}_\fc(M)$ if it decomposes 
into terms all with fiber degree $k$.
$D_M$ acts on the pair $(\alpha, \beta)$ as the matrix of operators
\begin{equation}
\left(
\begin{array}{cc} 
0 & -x \del_x - (f-k)  \\
x \del_x + k  & 0
\end{array} 
\right)
+ 
\left(\begin{array}{cc}
D_B + x^{-1} D_F & 0 \\ 0 & - D_B - x^{-1} D_F \end{array} \right).
\label{eq:DYprodfc}
\end{equation}

As before, this splits as $D_M = \Pi_0D_M \Pi_0 \oplus \Pi_\perp D_M
\Pi_{\perp}$.  If $\alpha$ and $\beta$ are $(j,k)$- and $(j-1,k)$-forms,
respectively, then 
\[
\Pi_\perp D_M \Pi_{\perp}= x^{-1} \widetilde{D}, \qquad
\mbox{where}\qquad 
\widetilde{D}=
D_{N,sc} + 
\left(\begin{array}{cc}
\Pi_\perp D_F \Pi_{\perp} & (b-j+1-f+k) x \\ (k-j) x & - \Pi_\perp
D_F \Pi_{\perp} 
\end{array} \right).
\]
The diagonal terms in this final matrix are constant in $x$ and invertible
on
$\Pi_\perp  L^2\Omega^*_\fc$, and reasoning as before, for any $a \in
\RR$, the mapping 
\[
\Pi_\perp D_M \Pi_{\perp}: x^a \Pi_{\perp} L^2 \Omega_{\fc}^*(M)
\longrightarrow x^{a-1} \Pi_{\perp} L^2 \Omega_{\fc}^*(M)
\]
has bounded inverse, $G^a_\perp$. 

On the other hand, $\Pi_0 D_M \Pi_0 \in \mbox{Diff}^1_b(N;
\Omega^{*,k}_\fc\calH^*(F))$ is a
$b$-operator (it has no $x^{-1}d_F$ or $x^{-1}\delta_F$ terms!), and
equals
\begin{equation}
\left(
\begin{array}{cc} 
D_B & -x \del_x - (f-k)  \\
x \del_x + k  & -D_B
\end{array} 
\right) .
\label{eq:sysfc}
\end{equation}
This operator preserves fiber degrees, so we can reduce to any
fixed $\Omega^{*,k}_\fc(M)$, for example when computing
indicial roots. We have
\[
I_{(\Pi_0 D_M \Pi_0)^2}(\gamma) =
\left(
\begin{array}{cc} 
D_B^2  -(\gamma+f-k)(\gamma+k) & 0 \\
0  & D_B^2-(\gamma+f-k)(\gamma+k)
\end{array} 
\right) .
\]
The critical exponent in the surjectivity calculation is $\gamma = -f/2$,
and inserting this into the expression above gives $D_B^2+(k-f/2)^2$
in both diagonal components. Hence elements in the nullspace
are in $L^2\calH^*(B)$ and are of fiber degree $k=f/2$. As before,
at this point one also checks that $-f/2$ is 
not an indicial root of multiplicity 
two, which simply involves showing as 
usual that (\ref{eq:sysfc}) has no solutions
of the form $\omega' x^{-f/2}\log x + \omega'' x^{-f/2}$.
 
In any case, so long as $a$ is not in the indicial 
set of $\Pi_0 D_M \Pi_0$, then
\[
\Pi_0 D_M \Pi_0: x^a \Pi_0 L^2 \Omega^*_\fc(M) \longrightarrow x^a
\Pi_0 L^2\Omega^*_\fc(M)
\]
is Fredholm, with generalized inverse $G^a_0$.

Altogether, this gives the generalized inverse
\[
G^a = G^a_0 \oplus G^a_\perp: 
x^{a}\Pi_0 L^2\Omega^*_{\fc}(M) \oplus
x^{a-1}\Pi_{\perp}L^2\Omega^*_{\fc}(M)
\longrightarrow x^{a} H^1\Omega^*_{\fc}(M),
\] 
and $I-G^a D_M = \Pi_M^a$ is the projection onto the nullspace of
$\Pi_0 D_M \Pi_0$ at weight $a$, all elements of which are
polyhomogeneous. 

In summary, the mappings
\begin{equation}
D_M: x^a H^1\Omega^*_{\fc}(M) \longrightarrow
x^a \Pi_0 L^2\Omega^*_{\fc}(M) \oplus
x^{a-1}\Pi_\perp L^2\Omega^*_{\fc}(M)
\label{eq:prfc1}
\end{equation}
and
\begin{equation}
D_M: x^a \Pi_0 H^1\Omega^*_{\fc}(M) \oplus x^{a+1} \Pi_{\perp} H^1
\Omega^*_{\fc}(M) \longrightarrow
x^a L^2\Omega^*_{\fc}(M)
\label{eq:prfc2}
\end{equation}
are Fredholm when $a$ is not an indicial root of $D_N$.

As in the fibred boundary case, if $\eta = D_M \zeta$
where $\eta \in \calA^a \Omega^*_\fc(M)$ and $\zeta \in
x^{c}\Pi_0 H^1\Omega^*_\fc(M) \oplus x^{c+1}\Pi_\perp H^1_\fc(M)$ 
for some $c<a$, then $\zeta \in \Pi_0 \calA^*_{\phg}\Omega^*_{\fc}(M) 
+ \calA^{a} \Omega^*_{\fc}(M)$. 

\subsection{Manifolds with nonproduct fiber bundle ends}
 
\subsubsection{Geometry of fibrations}

In this section we review some of the geometry associated to
a Riemannian fibration and use it to describe the precise structure of
$D_Y$. The exposition here is drawn from \S 10.1 of \cite{BGV}, \cite{dai},
\cite{Va}, and \cite{BC}, but since the notation in these sources varies
considerably, it has seemed worthwhile to develop this material in
detail. 

Suppose that $G = \phi^*(h) + k$ is a metric on the total space of a
fibration $Y$, where $\phi:Y \to B$ and $\phi^{-1}(b) = F_b$. As
before, we assume that $k$ annihilates the horizontal subbundle $T^H Y$,
which is the orthogonal complement of the vertical tangent bundle $T^V$,
and we let $P^V: TY \to T^V Y$, $P^H: TY \to T^H Y$ denote the orthogonal
projections. The tangent bundle $TB$ is naturally identified via $\phi_*$
with $T^H Y$, and we denote the lift of a section $X \in
\calC^\infty(B;TB)$ by $\tilde{X}$. In the following, we denote
sections of $T^V Y$ and $T^H Y$ by $U_1, U_2, \ldots$, and
$\tilde{X}_1, \tilde{X}_2, \ldots$, respectively. Finally, let
$\nabla^L$ denote the Levi-Civita connection of $G$.

The extent to which these subbundles fail to be parallel with respect
to $\nabla^L$ is measured in terms of two tensor fields, the second
fundamental form of the fibres and the curvature of the horizontal
distribution. The second fundamental form is the symmetric bilinear
form on $T^V Y$ defined by
\begin{equation}
\II_{\tilde{X}}(U_1,U_2) =
\big\langle \nabla_{U_1}^L U_2, \tilde{X} \big \rangle.
\label{eq:defII}
\end{equation}
We let $\II(U_1,U_2)$ be the horizontal vector given by
$$
\langle \II(U_1,U_2), \tilde{X} \rangle = \II_{\tilde{X}}(U_1,U_2),
$$
and we let $\II_{\tilde{X}}(U_1)$ denote the vertical vector
determined by 
$$
\langle \II_{\tilde{X}}(U_1),U_2 \rangle = \II_{\tilde{X}}(U_1,U_2).
$$
The nonintegrability of the horizontal distribution is measured by its
curvature, 
\begin{equation}
\calR(\tilde{X}_1,\tilde{X}_2) = P^V([\tilde{X}_1,\tilde{X}_2]),
\label{eq:curv1}
\end{equation}
which is tensorial and vertical. We define the horizontal vector
$\hat{\calR}_U(\tilde{X}_1)$ by
\begin{equation}
\big\langle \hat{\calR}_U(\tilde{X}_1),\tilde{X}_2\big\rangle
= \big\langle \calR(\tilde{X}_1,\tilde{X}_2),U\big\rangle
= \big\langle [\tilde{X}_1,\tilde{X}_2],U\big\rangle.
\label{eq:curv2}
\end{equation}
Four additional facts are used repeatedly. First, the bracket of
a vertical vector field with the horizontal lift of a vector field
from $B$ is again vertical, i.e.
\[
[\tilde{X},U] \in \calC^\infty (Y, T^V Y).
\]
This is proved by noting that vertical vector fields
are characterized by the fact that they annihilate functions
of the form $\phi^* f$, $f \in \calC^\infty(B)$.
Second, the Koszul formula determines the Levi-Civita
connection in terms of the metric and Lie brackets:
\[
\big\langle \nabla^L_{V_1}V_2,V_3 \big\rangle =
\frac12\left\{\big\langle [V_1, V_2], V_3 \big\rangle - \big\langle [V_2,
V_3] , 
V_1 \big\rangle + \big\langle [V_3,V_1], V_2\big\rangle  + \right.
\]
\[
\left. V_1\big\langle V_2,V_3\big\rangle + V_2 \big\langle
V_1,V_3\big\rangle -
V_3 \big\langle V_1,V_2 \big\rangle \right\},
\]
for any $V_1,V_2,V_3 \in \calC^\infty(Y,TY)$. 

Third, by definition of the induced Levi-Civita connection
$\nabla^F$ on the fibres,
\[
\big\langle \nabla^L_{U_1}U_2, U_3\big\rangle = \big\langle
\nabla^F_{U_1}U_2,
U_3 \big\rangle.
\]
Finally, since vertical and horizontal vector fields are perpendicular and
because the vertical distribution is integrable,
\[
\big\langle[U_1,U_2],\tilde{X}\big\rangle =
U_1\big\langle\tilde{X},U_2\big\rangle =
U_2\big\langle\tilde{X},U_1\big\rangle = 0.
\]
We now determine the vertical and horizontal components of
$\nabla^L_{V_1}V_2$, when the $V_j$ are, successively, vertical
and horizontal fields.
First, by definition, the horizontal part of $\nabla^L_{U_1}U_2$ is
\[
\big\langle \nabla^L_{U_1}U_2,\tilde{X}\big\rangle =
\big\langle \II(U_1, U_2), \tilde{X} \big\rangle,
\]
and also by definition, the vertical part is $\nabla^F_{U_1}U_2$.

{}From the Koszul formula and the expansion of
$\tilde{X}\big\langle U_1,U_2 \big\rangle$ using (\ref{eq:defII}),
\[
\big\langle \nabla^L_{\tilde{X}}U_1, U_2 \big\rangle =
\big\langle [\tilde{X},U_1],U_2\big\rangle -
\big\langle\II_{\tilde{X}}(U_1),
U_2\big\rangle,
\]
or in other words,
\[
P^V \nabla^L_{\tilde{X}}U = [\tilde{X},U] - \II_{\tilde{X}}(U).
\]
As for the horizontal component of $\nabla^L_{\tilde{X}}U$,
most of the terms in the Koszul formula vanish, leaving only
\[
\big\langle \nabla^L_{\tilde{X}_1}U,\tilde{X}_2\big\rangle =
-\frac12\big\langle [\tilde{X}_1,\tilde{X}_2],U\big\rangle =
-\frac12 \big\langle \hat\calR_U(\tilde{X}_1),\tilde{X}_2 \big\rangle.
\]

Next,
\[
\big\langle \nabla^L_{U_1}\tilde{X},U_2\big\rangle = -\big\langle
\tilde{X},\nabla^L_{U_1}U_2\big\rangle = -\big\langle \II_{\tilde{X}}
(U_1),U_2\big\rangle,
\]
is the vertical part of $\nabla^L_{U_1}\tilde{X}$ and the horizontal part is
\[
\big\langle \nabla^L_{U}\tilde{X}_1, \tilde{X}_2 \big\rangle
= \big\langle \nabla^L_{\tilde{X}_1}U +
[U,\tilde{X}_1],\tilde{X}_2\big\rangle
= \big\langle \nabla^L_{\tilde{X}_1}U, \tilde{X}_2 \big\rangle =
-\frac12 \big\langle \hat\calR_U(\tilde{X}_1),\tilde{X}_2 \big\rangle.
\]
where the Koszul formula is used for the final equality.

Finally, putting the covariant derivative on the other side of the
inner product and using the last equality of the previous displayed
formula, 
\[
\big\langle \nabla^L_{\tilde{X}_1}\tilde{X}_2,U\big\rangle =
\frac12 \big\langle \calR(\tilde{X}_1,\tilde{X}_2),U\big\rangle,
\]
and at last,
\[
\big\langle \nabla^L_{\tilde{X}_1}\tilde{X}_2,\tilde{X}_3\big\rangle
=
\big\langle \nabla^B_{X_1}X_2,X_3\big\rangle,
\]
where $\nabla^B$ is the Levi-Civita connection on $(B,h)$. This
last formula holds because all the terms in the Koszul formula
expansion only depend on $h$.

In summary, we have proved
\begin{proposition}
The Levi-Civita connection decomposes into vertical and  horizontal
components as
\begin{equation}
\begin{aligned}
\nabla^L_{U_1}U_2 & = & \nabla^F_{U_1}U_2 \quad & + & \quad\II(U_1,U_2) \\
\nabla^L_{\tilde{X}}U  & = & \ \ \left( [\tilde{X},U] -
\II_{\tilde{X}}(U)\right) \quad &-& \quad \frac12 \widehat\calR_U(\tilde{X})
\\
\nabla^L_{U}\tilde{X} & = & - \II_{\tilde{X}}(U) \quad &-& \quad
\frac12
\widehat\calR_U(\tilde{X}) \\
\nabla^L_{\tilde{X}_1}\tilde{X}_2 & = & \frac12 \calR(\tilde{X}_1,
\tilde{X}_2) \quad &+& \quad \left( \nabla^B_{X_1}X_2\right)\,\widetilde{}.
\end{aligned}
\label{eq:lcc}
\end{equation}
\end{proposition} 

We wish to define a new connection which preserves the splitting of $TY$. 
As a first guess, one might do this by projecting $\nabla^L$ onto the 
vertical and horizontal  
subspaces, i.e.\  to define $\nabla_{V}U = P^V (\nabla^L_VU)$, 
$\nabla_V(\tilde{X})= 
P^H(\nabla^L_V(\tilde{X})$, where $V$ is any vector (either
horizontal or vertical). The formul\ae\ above indicate which
terms should be subtracted from $\nabla^L$ to accomplish
this. However, there is another natural choice which turns out
to be more convenient for many computational purposes, given by
using the projected connection on the vertical bundle and
lifting the connection on the horizontal bundle from the Levi-Civita
connection on $B$. In other words we define
\[
\nabla := \left(P^V \nabla^L\right) \oplus \nabla^B,
\]
or even more specifically,
\[
\begin{array}{rcl}
\nabla_{U_1}U_2 & = & P^V (\nabla^L_{U_1}U_2), \\
\nabla_U \tilde{X} & = & 0,
\end{array}
\qquad
\begin{array}{rcl}
\nabla_{\tilde{X}}U &= &P^V(\nabla^L_{\tilde{X}}U) =
[\tilde{X},U] - \II_{\tilde{X}}(U), \\
\nabla_{\tilde{X}_1}\tilde{X}_2 & = & (\nabla^B_{X_1}X_2)~.
\end{array}
\]
We use this connection henceforth. Notice that it differs from
the projected connection only in the removal of the terms
$\frac12\widehat{\calR}_U(\tilde{X})$. One important feature of 
$\nabla$ vis a vis computations  
related to the families index theorem is that it is in `upper 
triangular form' with respect to the vertical/horizontal splitting, 
cf.\ \cite{BGV}. 

The difference tensor $Q = \nabla^L - \nabla$ is given by
\[
\begin{array}{rclrcl}
Q_{U_1}(U_2) &=& \II(U_1,U_2), \qquad &\qquad Q_{\tilde{X}}(U)
&=& -\frac12 \widehat\calR_U(\tilde{X}) \\
Q_{U}(\tilde{X}) &=& - \II_{\tilde{X}}(U) -\frac12
\widehat\calR_U(\tilde{X}) \qquad &\qquad
Q_{\tilde{X}_1}(\tilde{X}_2) & = & \frac12 \calR(\tilde{X}_1,
\tilde{X}_2).
\end{array}
\]
We note also that the torsion tensor of $\nabla$ is the negative of 
the skew-symmetrization of $Q$. 

We now express the deRham differential $d_Y$ and its adjoint in terms
of $\nabla$, $\II$ and $\calR$. Because the connections $\nabla^L$ and
$\nabla$ are both metric connections, they act on 1-forms by duality.  That
is, if $\phi$ is the 1-form given by $<w,\cdot>$, then $\nabla_Z \phi$ is
the
1-form given by $<\nabla_Z w, \cdot>$.  The action extends to forms of
higher
degree as a derivation.

Let $e_i$, $i=1, \ldots, f$ and $\eta_\mu$, $\mu = 1, \ldots, b$ be
orthonormal frame fields for $F$ and $B$, respectively, and
$\{e^i\}$, $\{\eta^\mu\}$ the dual coframe fields. It is standard that
\begin{equation}
d_Y = \sum_{i=1}^f e^i \wedge \nabla^L_{e_i} + \sum_{\mu=1}^b
\eta^\mu \wedge \nabla^L_{\eta_\mu},
\label{eq:defd}
\end{equation}
with analogous formul\ae\ for $d_F$ and $d_B$.
Now substitute $\nabla_L = \nabla + Q$ into (\ref{eq:defd}) to get first
\begin{equation}
d_Y e^j = d_F e^j + \sum \eta^\mu \wedge \nabla_{\eta_\mu}e^j
- \sum \left( \big\langle \II_{\eta_\mu}(e_i),e_j\big\rangle
\eta^\mu \wedge e^i + \frac12  \big\langle
\calR(\eta_\mu,\eta_\nu),e_j\big\rangle
\eta^\mu \wedge \eta^\nu\right),
\label{eq:dyej}
\end{equation}
and then
\begin{equation}
d_Y \eta^\mu = d_B \eta^\mu.
\label{eq:dyetamu}
\end{equation}
The last formula initially has many terms, all of which cancel,
but the result is no surprise since $d_Y \phi^* = \phi^* d_B$.

Now extend to forms of higher degrees. First, the splitting of $TY$
induces a decomposition
\[
\Lambda^k(T^* Y) = \bigoplus_{p+q=k}\Lambda^{p,q}(T^*Y),
\qquad
\mbox{where}\qquad
\Lambda^{p,q}(T^*Y) = \Lambda^p((T^VY)^*) \otimes
\Lambda^q((T^HY)^*).
\]
We regard the space of sections $\Omega^{p,q}(Y)$ as 
the completed tensor product $\Omega^p(B)\, \hat\otimes \, \Omega^q(Y,T^VY)$.
By construction, $\nabla$ preserves this splitting.
Thus for $\omega \in \Omega^{p,q}(Y)$, with $\omega = \phi^*(\al)
\wedge \be$, $\al \in \Omega^p(B)$ and $\be \in \calC^\infty(Y,
\Lambda^q((T^VY)^*)$,
\[
d_F \left(\phi^*(\al) \wedge \be\right) = (-1)^p\phi^*(\al ) \wedge d_F \be
\]
and we also define
\[
\tilde{d}_B \phi^*(\al) \wedge \be =
\phi^*(d_B \al) \wedge \be + (-1)^p\phi^*(\al)\wedge
\left(\sum_{\mu}\eta^\mu \wedge \nabla_{\eta_\mu} \be \right).
\]
Rewrite (\ref{eq:dyej}) as
\[
d_Y e^j = d_F e^j + \tilde{d}_B e^j - \II(e^j) -\frac12 \calR(e^j),
\]
where
\[
\II(e^j) =  \II_{\mu i j}\, \eta^\mu \wedge e^i, \qquad
\calR(e^j) = \calR_{\mu\nu j}\eta^\mu \wedge \eta^\nu.
\]
To simplify notation, let ${\rm R} = -\frac12 \calR$.
Then we have proved the first part of the
\begin{proposition} $d_Y = d_F + \tilde{d}_B - \II + {\rm R}$,
\ \ $\delta_Y = \delta_F + (\tilde{d}_B)^* - \II^* + {\rm R}^*$.
\label{pr:dYdY}
\end{proposition}
The second part is tautologous. Notice that
\[
\begin{array}{rcl}
d_F: \Omega^{p,q}(Y) \to \Omega^{p,q+1}(Y), & \qquad &
\tilde{d}_B: \Omega^{p,q}(Y) \to \Omega^{p+1,q}(Y)\\
\II: \Omega^{p,q}(Y) \to \Omega^{p+1,q}(Y), &\qquad&
{\rm R}: \Omega^{p,q}(Y) \to \Omega^{p+2,q-1}(Y).
\end{array}
\]

We can deduce some useful information from the fact that both $d_Y$ and 
$d_F$ are legitimate differentials, i.e.\ their squares are zero.
First, there is a Kodaira decomposition on the fibres, so any smooth
form $\alpha$ on $Y$ can be decomposed uniquely and orthogonally as
$\alpha = d_F \eta + \delta_F\mu + \gamma$, where $\gamma$ is fibre
harmonic.  Thus 
\begin{equation}
\Pi_0 d_F = \Pi_0 \delta_F = d_F \Pi_0 = \delta_F \Pi_0 =0.
\label{eq:reldf}
\end{equation}
Second, applying $d_Y^2=0$ to a form of pure bidegree and
decomposing into bidegrees gives
\begin{equation}
\begin{array}{rcl}
{\rm R}^2 & = & 0 \\
d_F^2 & = & 0 
\end{array}
\qquad
\begin{array}{rclrcl}
d_F(\tilde{d}_B - \II) & + &  (\tilde{d}_B - \II)d_F  & = &  0 \\
{\rm R} (\tilde{d}_B - \II) & + & (\tilde{d}_B - \II) {\rm R} &=& 0\\
d_F {\rm R} & + & {\rm R} d_F & = & - (\tilde{d}_B - \II)^2,
\end{array}
\label{eq:rels}
\end{equation}
with analogous relationships between the adjoints of these operators.

Now define the operator
\[
\frakd = \Pi_0\big(\tilde{d}_B - \II\big)\Pi_0;
\]
this acts on the space of fibre-harmonic forms.
\begin{proposition} The operator $\frakd$ and its adjoint $\frakd^*$
are differentials, i.e.\ $\frakd^2 = (\frakd^*)^2 = 0$.
\label{pr:frakddiff}
\end{proposition}
\begin{proof}
It suffices to prove only one of these.
Recalling that $\Pi_0 = I - \Pi_\perp$,
we have
\[
\frakd^2 = \Pi_0 \big( \tilde{d}_B - \II\big)^2\Pi_0 -
\Pi_0\big(\tilde{d}_B - \II\big) \Pi_{\perp} \big(\tilde{d}_B - \II\big)
\Pi_0;
\]
substituting from (\ref{eq:rels}) and using (\ref{eq:reldf}), this equals
\[
-\Pi_0 \big( d_F {\rm R} + {\rm R}d_F + (\tilde{d}_B - \II)
\Pi_{\perp} (\tilde{d}_B - \II) \big) \Pi_0
=- \Pi_0 \big(\tilde{d}_B-\II\big)\Pi_{\perp}\big(\tilde{d}_B -\II\big)
\Pi_0.
\]
Now, $d_F\big(\tilde{d}_B - \II\big)\Pi_0 = -\big(\tilde{d}_B - \II\big) d_F
\Pi_0 = 0$, so for any form $\alpha$, $(\tilde{d}_B - \II) \Pi_0
\alpha = d_F\eta + \gamma$,
with $\gamma$ fiber harmonic, and hence $\Pi_{\perp}\big(\tilde{d}_B -
\II\big) \Pi_0\alpha=d_F\eta$.  Finally,
\[
\frakd^2 \alpha = -\Pi_0 \big(\tilde{d}_B-\II\big) d_F \eta =
\Pi_0 d_F \big(\tilde{d}_B - \II\big) \eta =0.
\]
\end{proof}

\begin{corollary} Let $\DD = \frakd + \frakd^*$, and suppose that
$\DD \, \alpha = 0$ for some fibre-harmonic form $\alpha$. Then
$\frakd \alpha = \frakd^* \alpha = 0$, and so the terms $\alpha_{p,q}$
of pure bidegree also satisfy $\DD \alpha_{p,q} = 0$.
\end{corollary}
This follows just as for the usual Hodge Laplacian, for $\DD^2 =
\frakd^* \frakd + \frakd \frakd^*$ preserves bidegree and so
\[
0 = \langle \DD^2 \alpha ,\alpha \rangle = ||\frakd \alpha||^2
+ ||\frakd^* \alpha ||^2;
\]
in addition, we have that both $\frakd$ and $\frakd*$ commute 
with $\DD$.

\subsubsection{Hodge-de Rham operators in general}
The structure of the Hodge-de Rham operators for general exact fibred boundary
and fibred cusp metrics is obtained by substituting the expression for $d_Y$
from Proposition \ref{pr:dYdY} into (\ref{eq:DYprodfb}) and
(\ref{eq:DYprodfc}). To distinguish them, we write $D_\fb$ for
the operator $D_M$ associated to the fibred boundary metric $g_\fb$ and
$D_\fc$ for this operator associated to the fibred cusp metric $g_\fc$.
The action of $D_\fb$ on $\omega = \alpha/x^k + (dx/x^2)\wedge
\beta/x^{k-1}\in \Omega^{k,*}_\fb(M)$ is given by replacing the
second matrix in (\ref{eq:DYprodfb}) with
\[
\left(
\begin{array}{cc} 
D_F + x D_B - x(II + II^*) + x^2({\rm R} + {\rm R}^*) & 0
\\ 0 & -D_F - x D_B + x(II + II^*) - x^2({\rm R} + {\rm R}^*)
\end{array} 
\right).
\]
Similarly, the action of $D_\fc$ on $\omega = x^k \alpha + (dx/x)
\wedge x^k \beta \in \Omega^{*,k}_\fc(M)$ is obtained by substituting
\[
\left(
\begin{array}{cc} 
x^{-1} D_F + D_B - (II + II^*) + x({\rm R} + {\rm R}^*) & 0
\\ 0 & -x^{-1}D_F -  D_B + (II + II^*) - x({\rm R} + {\rm R}^*)
\end{array} 
\right) 
\]
for the second matrix in (\ref{eq:DYprodfc}).

As explained in the beginning of this section, the construction of
parametrices for $D_\fb$ and $D_\fc$ requires the machinery of
fibred-boundary pseudodifferential operators. The basic strategy is the
same in that one inverts $\Pi_0 D \Pi_0$ and $\Pi_\perp D \Pi_\perp$
separately, but now must also show that the off-diagonal
terms $\Pi_0 D \Pi_\perp$ and $\Pi_\perp D \Pi_0$, which no longer vanish,
play only an insignificant role. This is all carried out by Vaillant
\cite{Va},
cf.\ especially Proposition 3.27 there (although beware that the
Fredholm result is misstated in the special case $\lambda_0 = 0$) and
we shall simply quote the two results
we need, looking back to the product case for motivation.
Before stating these we remark that the operators $\Pi_0$,
$\Pi_\perp$ are only defined right at the boundary. However, the
fibred boundary structure requires that we have fixed the one-jet
of a definining function $x$ along the fibres, and this implies
that the spaces $x^c \Pi_ 0 L^2 \oplus x^{c \pm 1} \Pi_\perp L^2$
are well-defined for any $c \in \RR$ (because the weights only differ by $1$).

\begin{proposition} Suppose that $a$ is not an indicial root for
$\Pi_0 D_\fb \Pi_0$. Then
\begin{equation}
D_\fb: x^a H^1_\fb(M) \longrightarrow x^{a+1}\Pi_0 L^2\Omega^*_\fb(M)
\oplus x^{a}\Pi_\perp L^2\Omega^*_\fb(M)
\end{equation}
and
\begin{equation}
D_\fb: x^{a-1}\Pi_0 H^1\Omega^*_\fb(M) \oplus x^a \Pi_\perp H^1_\fb(M)
\longrightarrow x^a L^2\Omega^*_\fb(M)
\label{eq:gen1}
\end{equation}
are Fredholm. If $D_\fb \omega = 0$, then $\omega$ is polyhomogeneous,
with exponents in its expansion determined by the indicial roots
of $\Pi_0 x^{-1}D_\fb \Pi_0$, while if $\eta \in \calA^a \Omega^*_\fb(M)$,
$\zeta \in x^{c-1}\Pi_0 H^1\Omega^*_\fb(M) \oplus x^{c}
\Pi_\perp H^1_\fb(M)$ for $c<a$ and $\eta = D_\fb \zeta$, then
$\zeta \in \Pi_0 \calA^I_{phg} \Omega^*_{fb}(M) + \calA^{a}
\Omega^*_{fb}(M)$. 
\label{pr:fbmphg}
\end{proposition} 

\begin{proposition} Suppose that $a$ is not an indicial root for
$\Pi_0 D_\fc \Pi_0$. Then
\begin{equation}
D_\fc: x^a H^1_\fc(M) \longrightarrow x^{a}\Pi_0 L^2\Omega^*_\fc(M)
\oplus x^{a-1}\Pi_\perp L^2\Omega^*_\fc(M)
\end{equation}
is Fredholm. If $a+1$ is not an indicial root, then
\begin{equation}
D_\fc: x^a \Pi_0 H^1\Omega^*_{\fc}(M) \oplus x^{a+1} \Pi_{\perp} H^1
\Omega^*_{\fc}(M) \longrightarrow
x^a L^2\Omega^*_{\fc}(M)
\label{eq:gen2}
\end{equation}
is Fredholm.
If $D_\fc \omega = 0$, then
$\omega$ is polyhomogeneous, with exponents in its expansion determined by
the
indicial roots of $\Pi_0 D_\fc \Pi_0$, while if $\eta \in
\calA^a \Omega^*_\fc(M)$, $\zeta \in
x^{c}\Pi_0 H^1\Omega^*_\fc(M) \oplus x^{c+1}
\Pi_\perp H^1_\fc(M)$ where $c<a$ and $\eta = D_\fc \zeta$, then
$\zeta \in \Pi_0 \calA^I_{phg} \Omega^*_{fc}(M) + \calA^{a}
\Omega^*_{fc}(M)$. 
\label{pr:fcmphg}
\end{proposition} 

We remark that the indicial roots for the operators $\Pi_0 D \Pi_0$,
$D = D_\fb$ or $D_\fc$, are different than in the product case because
of the term $\II + \II^*$ and because of the action of $\tilde{d}_B$ on the
fiber part of forms; on the other hand the term
$\calR +
\calR^*$  is lower order at $x=0$ and does not affect the indicial roots.

\subsection{Hodge theorems for fibred boundary and fibred cusp metrics}

We now complete the proofs of the identifications of $L^2$ harmonic
forms with weighted cohomology in the two cases.

\medskip
\noindent{\bf Theorem 1C.} {\it If $(M,g)$ is a manifold with fibred 
boundary metric, 
then for every $k$ there is a natural isomorphism
\begin{equation}
L^2\calH^k(M) \longrightarrow \mbox{\rm Im}\,\big(\calW H^k(M,g_{\fb},\e)
\longrightarrow \calW H^k(M,g_{\fb},-\e)\big).
\label{eq:imsp}
\end{equation}
}
\medskip

\begin{proof}
If $\omega \in L^2\calH^k(M)$, then Proposition \ref{pr:fbmphg}
shows that $\omega$ is polyhomogeneous, and hence lies in
$x^{\e_0}L^2\Omega^k_\fb(M)$ for some $\e_0 > 0$ (with polyhomogeneous
coefficients). This gives the mapping
\[
L^2\calH^k(M) \longrightarrow \calW H^k(M,g_\fb,\e)
\longrightarrow 
\mbox{\rm Im}\,\big(\calW H^k(M,g_{\fb},\e)
\longrightarrow \calW H^k(M,g_{\fb},-\e)\big).
\]

If $[\omega] = 0$, then $\omega = d\zeta$ for some $\zeta \in
x^{-\e-1}L^2\Omega^{k-1}_\fb(M)$; by the discussion in \S 2.4,
we can choose $\zeta$ to be conormal. Write
\[
\omega = \sum_{p,q} \frac{\alpha_{p,q}}{x^p} + \frac{dx}{x^2} \wedge
\frac{\beta_{p,q}}{x^p}, \qquad
\zeta = \sum_{p,q} \frac{\mu_{p,q}}{x^p} + \frac{dx}{x^2} \wedge
\frac{\nu_{p,q}}{x^p},
\]
where $|\alpha_{p,q}|,|\beta_{p,q}|=\calO(x^{\frac{b+1}{2} + \epsilon_0})$
and $|\mu_{p,q}|,|\nu_{p,q}|=\calO(x^{\frac{b-1}{2} + \epsilon})$.
The usual integration by parts gives
\[
||\omega||^2 = \int_M d\zeta \wedge *\omega = \int_M d(\zeta \wedge *
\omega) = \lim_{x \rightarrow 0} \int_{B \times F} \zeta \wedge *\omega
= \lim_{x \rightarrow 0}\sum_{p,q} \int_{Y} \frac{\mu_{p,q}}{x^p}
\wedge \frac{*_Y \beta_{p,q}}{x^{b-p}},
\]
which vanishes, by the decay properties of the $\mu_{p,q}$ and
$\beta_{p,q}$. Thus $\omega = 0$, and this proves injectivity.

For surjectivity, we note that for sufficiently small $\epsilon >0$, the
space $L^2 \calH^*(M)$ can be identified with the cokernel of the map
\[
D_{\fb}: x^{-\epsilon-1} \Pi_0 H^1 \Omega^*_{\fb}(M) \oplus x^{-\epsilon}
\Pi_{\perp}
H^1\Omega^*_{\fb}(M) \longrightarrow x^{-\epsilon}L^2 \Omega^*_{\fb}(M).
\]
Thus we can write
\[
x^{-\epsilon} L^2 \Omega^*_{\fb}(M) = \mbox{\rm Im}
( D_{\fb}|_{x^{-\epsilon-1} \Pi_0 H^1 \Omega^*_{\fb}(M) + x^{-\epsilon}
\Pi_{\perp} H^1\Omega^*_{\fb}(M) }) \oplus L^2\calH^*(M).
\]
So suppose that
$\eta \in x^{\e}L^2\Omega^k_\fb(M)$ is a polyhomogeneous representative for a
class in the space on the right in (\ref{eq:imsp}). Then 
$\eta = D_\fb \zeta + \gamma$, where
$\zeta \in  \Pi_0 \calA^*_{\phg} \Omega^*_{\fb}(M) \oplus \Pi_\perp
\calA^{\epsilon}\Omega^*_{\fb}(M)$ 
and $\gamma \in L^2\calH^*(M)$.  In fact, comparing orders of vanishing
in $x$, we see that $\zeta = \zeta_0 + \zeta'$,
$\zeta' \in \calA^{\e}\Omega_{\fb}^*(M)$, and $\zeta_0 \in
\mbox{ker}\,I_{\Pi_0 x^{-1}D_M \Pi_0}((b-1)/2)$.  

We must analyze the structure of $\zeta_0$ more closely.
Acting on pairs $(\alpha, \beta)$, the indicial operator has the form
\[
I_{\Pi_0 x^{-1} D_\fb \Pi_0}((b-1)/2)= \left(
\begin{array}{cc} 
\DD & N_1 - (b-1)/2   \\
N_2 + (b-1)/2   & -\DD
\end{array} 
\right),
\]
where $\DD = \frakd + \frakd^*$. The operators $N_1$ and $N_2$ are defined
by $N_1 \beta_k = (b-k)\beta_k$ and $N_2 \alpha_k = -k \alpha_k$
(which agrees with the scattering case since $n=b+1$ there).
Following the calculation and reasoning for the scattering
case, we expand in terms of an eigenbasis for $\DD^2$ and deduce that
this indicial root has rank $1$ and that an element of the nullspace of
this indicial operator has the form
\[
\zeta_0 = x^{(b-1)/2}\left(\frac{\alpha_{(b-1)/2}}{x^{(b-1)/2}}
+ \frac{dx}{x^2} \wedge \frac{\beta_{(b+1)/2}}{x^{(b+1)/2}}\right)
\]
where $\alpha_{(b-1)/2},\ \beta_{(b+1)/2} \in \mbox{ker}\, \DD$.

We now have
\[
||\delta \zeta||^2 = \ <\eta - d \zeta - \gamma, \delta \zeta> \ =
\ <d(\eta - d \zeta - \gamma), \zeta> \ =
\lim_{x\to 0} \int_Y \zeta_0 \wedge d* \zeta_0
\]
\[
= \lim_{x\to 0} \int_Y \alpha_{(b-1)/2}\wedge d_Y *_Y \beta_{(b+1)/2}
= \ <\alpha_{(b-1)/2}, \frakd^* \beta_{(b+1)/2} +
{\rm R}^* \beta_{(b+1)/2}>_Y.
\]
But ${\rm R}^* \beta_{(b+1)/2}$ is a $((b-3)/2, *)$ form, so it pairs
trivially with $\alpha_{(b-1)/2}$, so this vanishes.

The rest of the argument is as in the scattering case.
\end{proof}

\medskip
\noindent{\bf Theorem 2C.} {\it If $(M,g_\fc)$ is a manifold with fibred 
cusp metric, then there is a natural isomorphism
\begin{equation}
L^2\calH^*(M) \longrightarrow \mbox{\rm Im}\,\big(WH^*(M,\e)
\longrightarrow WH^*(M,-\epsilon)\big).
\label{eq:imfcg}
\end{equation}
}

\begin{proof}
The proofs of the existence of this mapping and its injectivity are
nearly identical to those in the fibred boundary case, so we omit them.

For the surjectivity argument, we decompose
\[
x^{-\e}L^2\Omega^*(M) = \left( \left. \mbox{ran}\, D_M
\right|_{x^{-\e}\Pi_0 H^1_b \oplus x^{-\e+1}\Pi_{\perp}H^1_b}
\right) \oplus 
\left( \left. \mbox{ran}\, D_M \right|_{x^{-\e}\Pi_0 H^1_b \oplus
x^{-\e+1}\Pi_{\perp}H^1_b}
\right)^\perp.
\]
So we can write any $\eta \in x^{\e}L^2\Omega^*_\fc(M)$ which represents a
nontrivial class as $\eta = D \zeta + \gamma$,
where $\gamma \in L^2\calH^*(M)$ and $\zeta \in x^{-\e}L^2\Omega^*_\fc(M)$.
Since the indicial root $\gamma = -f/2$ occurs with multiplicity $1$, 
we have $\zeta = \zeta_0 + \zeta'$ where $\zeta' \in \calA^{\epsilon}
\Omega^*_{\fc}(M)$ and
\[
\zeta_0 = \sum_k (x^k\alpha_k + \frac{dx}{x} \wedge x^k \beta_k)x^{-f/2},
\]
where $\alpha_k$ and $\beta_k$ are independent of $x$ and $dx$. 
Matching up powers of $x$ in $\eta = D \zeta + \gamma$, we find that 
$\zeta_0$ is in the nullspace of the operator $I_{\frakd'}$, 
$\frakd' = \Pi_0 D_\fc \Pi_0$, which acts on $(*,k)$ forms by
\[
I_{\frakd'}(-f/2)= \left(
\begin{array}{cc} 
\DD & k - f/2  \\ k - f/2  & \DD
\end{array} 
\right).
\]
This implies that
$\alpha$ and
$\beta$ must both be forms on $B$ with coefficients in
$\mathcal{H}^{f/2}(F)$
and in the
kernel of $\DD$.  Thus the boundary term in the integration by parts
vanishes
as in the fibred boundary case.
\end{proof}

\subsection{From weighted cohomology to intersection cohomology}
To prove our main theorems, we must relate the weighted cohomology
groups appearing in the statements of Theorems 1C and 2C to intersection
cohomology groups. Most of the  work has already been done in \S 2.3, so it
remains only to reinterpret the answers.

The statement for fibred cusp metrics is slightly simpler, so 
consider that case first. We have proved that when $(M,g)$ is
a fibred cusp metric, then
\[
L^2\calH^*(M) \cong \mbox{\rm Im}\,\big(WH^*(M,\e)\longrightarrow
WH^*(M,-\epsilon)\big).
\]
Using Proposition \ref{pr:wcih}, this is equivalent to
\[
L^2\calH^*(M) \cong \mbox{\rm Im}(I\!H^*_{[\epsilon +(f/2)]}(X,\hB)
\longrightarrow I\!H^*_{[-\epsilon +(f/2)]}(X,\hB)),
\]
where $X$ is the compactification of $M$ defined in the introduction. The
two 
spaces on the right correspond to intersection cohomology with the
middle perversities
\[
\underline{\frakm}(f+1)=\left\{ \begin{array}{ll}
\frac{f-1}{2} & f \ \ \mbox{odd}\\ \frac{f}{2}   & f \ \ \mbox{even}
\end{array} \right.
\qquad
\overline{\frakm}(f+1)=\left\{ \begin{array}{ll}
\frac{f-1}{2} & f \ \ \mbox{odd}\\
\frac{f}{2}-1  & f \ \ \mbox{even}
\end{array} \right.,
\]
respectively. This proves the main

\medskip
\noindent{\bf Theorem 2.} {\it 
Suppose $(M,g)$ is a manifold with fibred cusp metric.
Then
\[
L^2\calH^*(M) \cong \mbox{\rm Im}\big(I\!H^*_{\underline{\frakm}}(X,B)
\longrightarrow I\!H^*_{\overline{\frakm}}(X,B)\big).
\]
}

We remark on a few special cases of this result:

\medskip

\noindent If $f=0$ (i.e.\ $(M,g)$ has cylindrical ends), then
\[
I\!H^*_{[\epsilon +(f/2)]}(X,\hB) = H^*(M,\del M), \qquad
I\!H^*_{[-\epsilon +(f/2)]}(X,\hB)) = H^*(M),
\]
and so we recover the image of relative in absolute cohomology,
as already proved in \S 4.

\medskip

\noindent If $\dim F = f > 0$, then the two spaces on the right
coincide when $f$ is odd, or even if we only have 
$H^{f/2}(F)=0$, i.e.\ $(X,B)$ is a Witt space. In either case, 
$L^2\calH^*(M)$ equals the (unique) middle perversity intersection 
cohomology $I\!H^*_{\frakm}(X,B)$.

We can see this simplification directly from the analysis in the last 
section. Recall the decomposition $\eta = d\zeta + \gamma$ for 
the closed form $\eta \in \calA^{\e}\Omega_\fc^k(M)$. 
We have $\zeta = \zeta_0 + \zeta'$ where $\zeta' \in \calA^\e
\Omega^*_\fc(M)$ and $\zeta_0$ is the sum of pullbacks of form on $B$ 
wedged with an element of $H^{f/2}(F)$. But the assumption that $X$ is a 
Witt space gives $\zeta_0=0$, and hence $[\eta] = [\gamma]$
already in $WH(M,g_\fc,\e)$. Thus in this case 
\[
WH(M,g_\fc,-\e) = WH(M,g_\fc,\e) = WH(M,g_\fc,0) = L^2\calH^*(M),
\]
and all these spaces are finite dimensional. This already follows from
\cite{Z}, Corollary 2.34. Finally, the discussion in \S 2.3 shows
how to interpret this in terms of intersection cohomology.

However, when $H^{f/2}(F) \neq 0$ the unweighted $L^2$ cohomology
is infinite dimensional, and the two middle perversity intersection
cohomologies are different. In this case, some sort of more elaborate 
analysis, as we have carried out in this paper, is needed.

We obtain the Hodge theorem for fibred boundary metrics by
a translation from the fibred cusp case. To do this, first rewrite 
the isomorphism 
\[
L^2\calH^k(M) \cong \mbox{\rm Im}\big(\calW H^k(M,g_{\fb},\e) \longrightarrow
\calW H^k(M,g_{\fb},-\e)\big)
\]
in terms of weighted $L^2$ cohomology for the associated fibred cusp metric
$g_\fc = x^2 g_\fb$. This gives
\[
L^2\calH^k(M) \cong \mbox{\rm Im}\,\big((WH^k(M,g_{\fc},n/2-k+\e)
\longrightarrow WH^k(M,g_{\fc},n/2-k-\e)\big),
\]
and hence by Proposition \ref{pr:wcih} we get

\medskip

\noindent{\bf Theorem 1.} {\it If $(M,g)$ is a fibred boundary metric, 
then
\[
L^2\calH^k(M) \cong \mbox{\rm Im}\, (I\!H^k_{[\frac{n+f}{2}-k+\e]}
(X,\hB)\longrightarrow I\!H^k_{[\frac{n+f}{2}-k-\e]}(X,\hB)).
\]
}

\medskip

We list the various cases:

\medskip

\noindent Suppose $b$ is even. Since $n = b+f+1$, this is the same as
$n+f$ is odd, and then the two groups are the same, so that
\[
L^2\calH^k(M) \ \cong \ I\!H^k_{f+\frac{b}{2}-k}(X,\hB)
\ \cong \ \left\{ \begin{array}{ll}
H^k(X, B) & k \leq \frac{b}{2} \\
I\!H^k_{f-1}(X,B) & k = \frac{b}{2} +1 \\
\vdots  & \\
I\!H^k_{0}(X,B) & k=n-\frac{b}{2} +1 \\
H^k(X \setminus B) & k\geq n-\frac{b}{2}
\end{array} \right. .
\]
Just as in the fibred cusp case, when $b$ is even, the form
$\zeta_0$ which arises in the surjectivity argument
must vanish since it lies in $\Omega^{(b \pm 1)/2,*} = \{0\}$ on
the boundary. Hence the map $\Phi$ is now surjective onto
$\calW H^*(M,g_\fb,\e)$. In this case the range of $D$ is
closed, and the theorem follows from the techniques of \cite{Z}.

\medskip

\noindent When $b$ is odd, 
\[
L^2\calH^k(M) \ \ \cong \ \ \mbox{\rm Im}(I\!H^k_{f+\frac{b+1}{2}-k}(X,\hB)
\longrightarrow I\!H^k_{f+\frac{b-1}{2}-k}(X,\hB))
\]
\[
\cong \qquad \left\{ \begin{array}{ll}
H^k(X,\hB) & k \leq \frac{b-1}{2} \\
\mbox{\rm Im}\big(H^k(X,\hB) \longrightarrow I\!H^k_{f-1}(X,\hB)\big) & k =
\frac{b-1}{2} +1 \\
\mbox{\rm Im}\big(I\!H^k_{f-1}(X,\hB) \longrightarrow I\!H^k_{f-2}(X,\hB)\big)
& k = \frac{b-1}{2} +2 \\
\vdots  & \\
\mbox{\rm Im}\big(I\!H^k_{1}(X,\hB) \longrightarrow I\!H^k_{0}(X,\hB)\big) & k
= n-\frac{b-1}{2} -2 \\
\mbox{\rm Im}\big(I\!H^k_{0}(X,\hB) \longrightarrow H^k(X\setminus \hB)) & k =
n-\frac{b-1}{2} -1 \\
H^k(X\setminus \hB) & k\geq n-\frac{b-1}{2}
\end{array} \right. .
\]

\medskip

\noindent Simpler corollaries of this theorem, for cases when $F$ is a
sphere and $X$  a smooth manifold, were stated in Corollary 1 in the
introduction.

\section{Relationship to other works}

We now briefly discuss some consequences of the Hodge theorems
proved here and their relationship with other work in the field.

\medskip

\noindent{\bf Carron's Hodge theorem for manifolds with flat ends:} In a 
recent paper \cite{Car}, Carron has calculated the Hodge cohomology
for manifolds with finitely many ends, on all of which it is assumed 
that the curvature tensor vanishes identically. He uses two main tools: 
a precise geometric structure theorem for flat ends \cite{ES}, and 
his theory of nonparabolicity at infinity in order to obtain 
new function spaces, which are extensions of $H^1_0\Omega^*(M)$ and on 
which the range of $D$ is closed. This work 
has substantial overlap with ours in the sense that many but not
all fibred boundary and fibred cusp metrics are nonparabolic at
infinity and satisfy the extra conditions implied by the flatness
hypothesis.

\medskip

\noindent{\bf The signature formula of Dai and Vaillant} 
As discussed in the introduction, an immediate corollary of 
Theorems 1 and 2 is that
\begin{equation}
{\mbox{sgn}}_{L^2}(M,g) = {\mbox{sgn}}\,
\mbox{\rm Im}\,\big( I\!H_{\underline{\frakm}}(X,B)
\longrightarrow I\!H_{\overline{\frakm}}(X,B)\big).
\label{eq:sgnih}
\end{equation}
This formula holds both for fibred boundary and fibred cusp metrics. 

On the other hand,  
there is an $L^2$ signature theorem for manifolds with fibred cusp ends proved 
by Dai \cite{dai} and generalized by Vaillant \cite{Va}:
\begin{equation}
\mbox{sgn}_{L^2}(M,g) = \mbox{sgn}\,
\mbox{\rm Im}\,\big( H^*(M, \del M) 
\longrightarrow H^*(M))\big) + \tau.
\label{eq:sgnra}
\end{equation}
The final term here is the $\tau$ invariant, originally defined by Dai,
which is a sum of signatures coming from the higher terms
in the Leray spectral sequence for the fibration of $\del M$. 
Combining these two signature theorems now identifies
$\tau = \tau(\del M)$ with the difference of the two algebraic
signatures in (\ref{eq:sgnih}) and (\ref{eq:sgnra}), see
(\ref{eq:difftau}) in the introduction.  The original
definition of $\tau$ involves algebraic signatures on the
higher terms (i.e. the $E_k$ terms, $k \geq 3$) of the Leray
spectral sequence of the fibration for $\del M$. It seems very
tempting to conjecture that the summands in this definition 
arise from signatures on the weighted cohomology for weights
$\pm a$, where $a$ varies from some small positive number
to one sufficiently large so that the weighted cohomologies $WH(M,g,\pm a)$
equal the relative and absolute cohomologies, respectively. 
There should be finitely many jumps in this deformation,
and the intermediate weighted cohomologies should correspond
to intersection cohomologies with perversities varying from
lower middle or upper middle to one of the extremes.  We shall
return to a precise exploration of these ideas elsewhere.

\medskip

\noindent{\bf Hitchin's Hodge theorem:} 
The next section contains an explanation of our Hodge and
signature theorems in several interesting examples. Most of those 
examples are hyperk\"ahler, and the Hodge cohomology of such manifolds
has been recently studied by Hitchin \cite{hitchin2}.  Amongst
his results is one particularly relevant to our paper:

\medskip
\noindent{\bf Theorem 3 (Hitchin):}\ {\it Let $M$ be a complete 
hyperk\"ahler manifold of real dimension $4k$ such that one of the 
K\"ahler forms $\omega_i$ satisfies $\omega_i=d\beta$, where 
$\beta$ has linear growth. Then any $L^2$ harmonic form on $M$ is of degree
$2k$ and is self-dual or antiself-dual provided that $k$ is even 
(respectively, odd).}

\medskip

This implies 
\begin{corollary} If $M$ is a hyperk\"ahler manifold as above, then
$\dim L^2\calH^*(M,g)= |{\mbox{\rm sgn}}_{L^2}(M,g)|$.
\end{corollary}

Hence for the class of hyperk\"ahler manifolds satisfying the hypothesis
of Hitchin's theorem (including most of the examples in the next section),
the Hodge cohomology can be computed from the $L^2$-signature index 
theorem of Dai and Vaillant. 

We obtain two consquences which follow from this result and the analysis 
developed for the proofs of our main theorems. The first gives an interesting 
topological obstruction to the existence of a fibred boundary or fibred cusp 
hyperk\"ahler metric satisfying the linear growth hypothesis of 
Theorem 3. 

\begin{corollary} If $M$ is a hyperk\"{a}hler manifold as in Theorem
3 which is either of fibred cusp or fibred boundary type, 
then the intersection form on $H^*(M,\del M)$ is semidefinite so that
$\mbox{\rm sgn}(M)$ is nonpositive if $k$ is odd and nonnegative if $k$ 
is even.
\end{corollary}

\begin{proof} To be definite, suppose $g$ is a fibred cusp metric.
We know by Theorem 3 above that the intersection form on 
$L^2{\calH}^{2k}(M,g)$ is semidefinite of the correct sign. 
On the other hand, the topological signature of a manifold with boundary
is by definition the index of the intersection form on the image of 
(middle degree) relative cohomology in absolute. Thus we must show
that this latter intersection form is also semidefinite. 

Suppose that $\eta$ and $\nu$ are smooth closed compactly supported
$2k$-forms which represent nontrivial classes in $\mbox{\rm Im}(H^{2k}
(M,\del M)\to H^{2k}(M))$. By Theorem 2, or rather its proof in \S 5,
we have $\eta = d\zeta + \gamma$, $\nu = d\xi + \rho$ where
$\gamma, \rho \in L^2\calH^{2k}(M)$; we also have that
$\zeta = \zeta_0 + \zeta'$,  where $\zeta' \in \calA^\e\Omega^{2k}_\fc(M)$
and $\zeta_0$ is polyhomogeneous with growth at 
just the critical value for square integrability and in addition is
fiber harmonic form and in the kernel of $\DD$. There is a similar
decomposition for $\xi$. 

We now compute that
\[
\int_M \eta \wedge \nu = \int_M (d\zeta + \gamma) \wedge (d\xi + \rho)
=\int_M d\zeta \wedge d\xi + \int_M d\zeta \wedge \rho
+ \int_M \gamma \wedge d\xi + \int_M \gamma \wedge \rho.
\]
Now integrate by parts in each of the first three terms on the right;
using the information in the last paragraph, the boundary terms all
vanish, and we are left with the equality of the pairing of 
$\eta$ and $\nu$ with the pairing of $\gamma$ and $\rho$, as desired.
\end{proof}

\begin{remark} This topological obstruction is investigated further
in \cite{hausel-swartz} for toric hyperk\"ahler varieties.
\end{remark}

The argument in the proof above also yield 
\begin{corollary} If $M$ has a hyperk\"ahler fibred boundary metric as 
above, then the $\tau$ invariant of $\del M$ is 
non-positive if $k$ is odd and non-negative if $k$ is even.
\end{corollary}

\section{Examples}
A mathematically interesting theme in contemporary research in string theory 
involves the use of duality to predict the dimensions of spaces 
of $L^2$ harmonic forms on various classes of noncompact manifolds.
Probably the most famous of these is the S-duality conjecture made
by Sen in \cite{sen1}, which predicts the dimension of the Hodge
cohomology on moduli spaces of monopoles on $\RR^3$; these moduli spaces
include the Atiyah-Hitchin manifold, the Taub-NUT space and its 
higher dimensional generalizations. A similar S-duality prediction in 
\cite{vafa-witten} concerns the Hodge cohomology of quiver varieties,
while \cite{hausel} contains a mathematical conjecture about the Hodge 
cohomology of moduli of Higgs bundles. Similar to Sen's conjecture, these
last predictions equate the Hodge cohomology of these moduli spaces
with the image of compactly supported cohomology in absolute cohomology.
We also mention the predictions about Hodge cohomology in \cite{sen2},
for multi-Taub-NUT spaces, and in \cite{gomis-etal}, for the $G_2$ space
constructed in that paper.

The justification of these predictions has been a key motivation for 
our work. In this final section we examine these conjectures
in light of the results of this paper. The point is that,
particularly in the low dimensional cases, the moduli spaces
in these conjectures carry natural fibred boundary metrics, and
hence our Theorem 1 can be applied. We discuss several examples 
where we can confirm the predictions, but notably, we also show
that the $L^2$ harmonic form predicted to exist on the ALF $G_2$ 
space of \cite{gomis-etal} does not in fact exist. This is labeled 
as a $U(1)$-puzzle in Section 6 of that paper, and awaits further 
explanation.

Many of the calculations below have been or could be done using techniques
already in the literature. For example, Hitchin \cite{hitchin2} has already 
settled Sen's S-duality conjecture for the Atiyah-Hitchin and
Taub-NUT manifolds. Likewise, the computations for all hyperk\"{a}hler 
ALE spaces follow from Theorem 3 above and the computation of Hodge
cohomology in the b-case, which was previously known, \cite{APS}, 
\cite{Me-aps}. For spaces with hyperk\"{a}hler metrics of 
fibred boundary type, the calculations follow from Theorem 3 again
and the signature formula (\ref{eq:sgnra}) of Dai and Vaillant.
We hope the advantages of our more unified approach to these problems
is apparent and that our results give new topological insight 
even in the previously understood cases. We shall state as a corollary 
those applications which we believe are new. 

\subsection{Gravitational Instantons}
\label{gravinst}
A gravitational instanton is by definition, \cite{hawking}, a 
$4$-dimensional complete hyperk\"ahler manifold. In all known, topologically
finite and non-compact 
examples, the metric is of fibred boundary type. These examples can be
separated into three classes: ALE (short for asymptotically locally
euclidean), where $F$ is a point; ALF (short for asymptotically locally flat),
where $F=S^1$; and ALG (by induction) where $F=S^1\times S^1$.

The space $L^2\calH^2(M)$ of $L^2$ harmonic $2$-forms for gravitational 
instantons is particularly interesting since it contains the curvatures 
of $\mbox{U}(1)$ Yang-Mills connections. Because of this, we shall 
also mention 
what is known about $\mbox{SU}(2)$ Yang-Mills connections on gravitational 
instantons and how these $\mbox{U}(1)$ Yang-Mills connections fit into 
that picture as subspaces of reducible connections.

\subsubsection{ALE gravitational instantons}
In his thesis, Kronheimer classified all ALE gravitational 
instantons, \cite{kronheimer1}, \cite{kronheimer2}. The underlying
manifolds in this classification are (diffeomorphic to) minimal 
resolutions of $\CC^2/\Gamma$, where $\Gamma$ is a finite subgroup
of $\mbox{SU}(2)$. These are of type $A_k$, $D_k$, $E_6$, $E_7$ or $E_8$.
Denoting the resolution of $\CC^2/\Gamma$ by $M_\Gamma$, 
the correspondence is given by the fact that the intersection form on 
$H^2_c(M_\Gamma)$ is isomorphic to the Cartan matrix of some 
simply laced Lie algebra of type ADE. Topologically, this means that 
$M_\Gamma$ retracts to a configuration of Lagrangian $2$-spheres forming
the corresponding Dynkin diagram. The intersection form 
gives a pairing $H_c^2(M_\Gamma) \times H^2(M_\Gamma) \to \ZZ$, 
and since the Cartan matrix defining the form is always negative definite, 
we see that the forgetful map 
$H_c^2(M_\Gamma)\to H^2(M_\Gamma)$ is an isomorphism. 

Now apply Theorem 1 to get the well-known result that $L^2\calH^k(M)$
is nontrivial only in degree $2$, and 
\[
L^2\calH^2(M_\Gamma,g_{{\mathrm{ALE}}})\cong H^2(M_\Gamma).
\]
In particular, if $k$ is the number of conjugacy classes in $\Gamma$, then
$\dim L^2\calH^2(M_\Gamma,g_{\mathrm{ALE}})=k-1$. 

A nice explicit construction of $k-1$ independent elements giving
a basis of $L^2\calH^2(M_\Gamma)$ in this case appears in \cite{gocho-nakajima}. 
The paper \cite{kronheimer-nakajima} combines this with \cite{kronheimer1}
to construct all finite energy $\mbox{U}(k)$ Yang-Mills instantons on 
$M_\Gamma$.

\subsubsection{ALF gravitational instantons}
There is no classification known for ALF gravitational instantons
parallel to that of Kronheimer for the ALE case. However, recently 
Cherkis and Kapustin \cite{cherkis-kapustin2} have conjectured 
a classification scheme: using a physics argument they predict
that all ALF instantons are of the types: $A_k$, $D_k$, so that $D_0$ stands 
for the Atiyah-Hitchin manifold. 

Consider first the $A_k$ (for $k\geq 1$) and $D_k$ (for $k\geq 4$) families. 
The underlying manifolds 
of these gravitational instantons are the same as in the ALE case, 
although the metrics are of course now ALF. Thus now $\Gamma$ 
is either a cyclic or dihedral subgroup of $\mbox{SU}(2)$ and
$M_\Gamma$ the minimal resolution of $ \CC^2/\Gamma$.
The $A_k$ family was constructed first in \cite{hawking} (see below for 
the details), while the $D_k$ family appears in 
\cite{cherkis-kapustin2} and \cite{cherkis-kapustin1}.

The following corollary confirms the prediction made in \cite{sen2} 
concerning the Hodge cohomology of ALF gravitational instantons in 
the $A_k$ case, but includes the $D_k$ case as well. 

\begin{corollary} Suppose $\Gamma\subset \mbox{SU}(2)$ is a finite cyclic or 
dihedral subgroup, and let $k$ be the number of conjugacy classes in
$\Gamma$. If $(M_\Gamma,g_{\mathrm{ALF}})$ is the associated
ALF gravitational instanton, then $\dim L^2\calH^2(M_\Gamma) = k$;
$L^2\calH^d(M_\Gamma)$ is trivial for $d \neq 2$. 
\end{corollary}
\begin{proof} 
 In both the $A_k$ and $D_k$ settings
$\overline{X_\Gamma}=X_\Gamma\cup S^2$. The Mayer-Vietoris
sequence gives that $H^*(\overline{X_\Gamma})\cong 
H^2(X_\Gamma)\oplus H^0(S^2)$. Therefore, by (\ref{sphereeven}),
$\dim L^2\calH^2(M,g_{\mathrm{ALF}})= \dim H^2(X_\Gamma) +1=k$.

Alternatively, apply Theorem 3 and (\ref{eq:sgnra}). One calculates
that the $\tau$ invariant of the fibration at infinity is $-1$, hence
$\mbox{sgn}_{L^2}(M_\Gamma,g_{\mathrm{ALF}})=\mbox{sgn}(M_\Gamma)-1=
-k$. The result follows by applying Theorem 3 again. 
\end{proof}

A consequence of this result is that for an ALF gravitational instanton 
$M_\Gamma$ there is, up to scaling, a unique $L^2$ harmonic form; this
form is exact but not, of course, in the range of $d$ on $L^2$. In the 
$A_k$ case, the metric and all $L^2$ harmonic $2$-forms are known explicitly. 
We now explain this in more detail and determine which $L^2$ harmonic form 
is exact.

The explicit construction of the ALF gravitational instantons of type 
$A_k$ uses the Gibbons-Hawking ansatz \cite{gibbons-hawking}:
\[
g_{\mathrm{ALF}}=V(dx_1^2+dx_2^2+dx_3^2)+V^{-1}(d\theta+\alpha)^2,
\]
where $\alpha$ is a $1$-form on $\RR^3$ such that $d\alpha=*dV$. 
There is a metric $g_{\mathrm{ALF}}^k$ of this type which lives on a 
four-manifold $M_k$ and admits an isometric circle action with $k$ 
fixed points. Away from these fixed 
points, $M_k$ fibers over $\RR^3\setminus \{p_1,\dots,p_k\}$ with
$S^1$ fibers, and it induces a degree $-1$ 
fibration around each $p_i\in \RR^3$. 
Here $(x_1,x_2,x_3)$ is the standard coordinate system on $\RR^3$ and
$\theta \in S^1$. Finally, 
\[
V=\sum^k_1 \frac{2m}{|x-p_i|} + 1, \qquad m > 0.
\]
These are called Gibbons-Hawking or multi-Taub-NUT metrics, 
and $g^{1}_{\mathrm{ALF}}$ is the famous Taub-NUT metric.

The paper \cite{ruback} explicitly describes the $k$-dimensional space 
$L^2\calH^2(M_k)$ as follows: 
\[
\Omega_i=d\xi_i, \qquad i=1,\dots,k,
\] 
where
\[
\xi_i=\alpha_i-\frac{V_i}{V}\left(d\theta + \alpha\right),
\qquad V_i=\frac{2m}{|x-p_i|},
\qquad \mbox{and}\quad d \alpha_i=*d V_i.
\]
This description is only local in the given coordinate chart, and 
indeed, $\xi_i$ extends globally only as a connection on a
$\mbox{U}(1)$ bundle. Its curvature $\Omega_i$ is globally defined.
There is one exception: the connection $\xi = \sum \xi_i =
\frac{1}{V}(d\theta+\alpha)-d\theta$ is gauge equivalent to 
$\frac{1}{V}(d\theta+\alpha)$, which extends globally as 
the metric dual of the Killing vector field
$\frac{\partial}{\partial \theta}$ from the circle action.
Its curvature is the $L^2$ harmonic $2$-form 
$d\left(\frac{1}{V}(d\theta+\alpha)\right)$.
For the Taub-NUT metric, i.e.\ when $k=1$, this $2$-form was 
discovered by Gibbons \cite{gibbons} and exhibited as support for Sen's
S-duality conjecture. (As already noted, Hitchin \cite{hitchin2} 
settled Sen's conjecture in this case by proving that there are no other
non-trivial $L^2$ harmonic forms.) 

Our result explains the topological origin of Gibbons' $L^2$ harmonic 
$2$-form. For although $M_1$ is diffeomorphic to $\RR^4$, its 
compactification (as an ALF space) is ${X}_1=\CC P^2$.
The non-trivial cohomology of $\CC P^2$ in degree $2$ is the topological 
source of Gibbons' $L^2$ harmonic $2$-form.

The other infinite family of ALF gravitational instantons, of type
$D_k$, was constructed in \cite{cherkis-kapustin1,cherkis-kapustin2} as
moduli spaces of certain singular $\mbox{SU}(2)$ monopoles on $\RR^3$. 
The metrics are defined using twistor theory, so are not as explicit as 
the Gibbons-Hawking metrics above. However, for $k\geq 4$, 
Theorem 1 again gives a 
$k$-dimensional space of $L^2$ harmonic $2$-forms, a $1$-dimensional
subspace of which is exact. It would be interesting to find these
harmonic forms explicitly.

We now come to the Atiyah-Hitchin manifold $M$ \cite{atiyah-hitchin}.
As explained in \cite{hitchin1}, the compactification of this
space is obtained by adding a copy of $\RR P^2$, and in fact 
$M\cup \RR P^2=S^4$. Hence (\ref{sphereeven}) shows that $L^2\calH^*(M) = 0$.
However, $\pi_1(M)$ is $\ZZ_2$, and the universal cover
$\tilde M$ has compactification $\tilde{M} \cup \RR P^2=\CC P^2$.
Therefore $L^2\calH^2(\tilde{M})$ is one-dimensional. This $2$-form was 
constructed by Sen in \cite{sen1}, and Hitchin \cite{hitchin2} proved its uniqueness. 
Our proof of Sen's conjecture, through (\ref{sphereeven}), explains 
the topological origin of this form, since it comes from 
the $1$ dimensional $H^2(\CC P^2)$.

In contrast with the ALE case, very little is known about Yang-Mills 
instantons on these ALF gravitational instantons (though, of course, 
the discussion above can be applied to understand the situation for
$\mbox{U}(1)$ Yang-Mills instantons). Recently new families of $\mbox{SU}(2)$
Yang-Mills instantons on multi-Taub-NUT spaces have been found, cf.\ 
\cite{etesi-hausel2}, \cite{etesi-hausel3}. In particular, 
\cite{etesi-hausel3} contains an intrinsic construction of the
$L^2$ harmonic forms $\Omega_i$ defined above as the curvatures of 
reducible $SU(2)$ Yang-Mills instantons.
 
We conclude this section with a final example, the well-known 
Euclidean Schwarzschild space $M$, which is a complete Ricci-flat 
$4$-manifold but not hyperk\"ahler. Its Hodge cohomology is calculated in 
\cite{etesi-hausel1} using techniques from \cite{hitchin2}, and
it is shown there that $L^2\calH^k(M) = 0$ when $k \neq 2$ and
$L^2\calH^2(M)$ is $2$-dimensional, with a $1$-dimensional subspace of 
(anti)-self-dual solutions. This is explained neatly by
(\ref{sphereeven}): namely $M$ is diffeomorphic to $\RR^2 \times S^2$,
and is ALF with $F=S^1=\partial(\RR^2)$, hence it compactifies as 
$X = S^2 \times S^2$. Applying (\ref{sphereeven}), we see that
the Hodge cohomology of $M$ is concentrated in degree $2$, and
\[
\dim L^2\calH^2(M)= \dim H^2(X)= \dim H^2(S^2\times S^2)=2.
\]
As explained in \cite{etesi-hausel1}, the self-dual $L^2$ harmonic 
$2$-forms on $M$ had already appeared in the physics literature in 
the disguise of $\mbox{SU}(2)$ Yang-Mills instantons \cite{charap-duff}.

\subsubsection{ALG gravitational instantons}
The ALG gravitational instantons are the most recent of these spaces
to be studied and examples have only recently been constructed
\cite{cherkis-kapustin3}; they arise as moduli spaces of periodic monopoles
on $\RR^2\times S^1$. In these examples the underlying manifold $M$ is an 
elliptic fibration of type $D_1$, $D_2$, $D_3$, $D_4$ or $E_6$, $E_7$, $E_8$,
cf.\ \cite{cherkis-kapustin3} for the precise meaning of this.
They all have a fibred boundary metric with $F=T^2$, and hence their
compactification $X = M \cup S^1$ is not a Witt space.
Theorem 1 gives 
\begin{corollary} Let $(M,g_{\mathrm{ALG}})$ be an $ALG$ gravitational 
instanton. Then 
\[
L^2\calH^2(M,g_{\mathrm{ALG}}) \cong 
\mbox{\rm Im}\,(H^2(M,\del M)\to H^2(M)),
\]
is an isomorphism, or in other words, $\dim L^2\calH^2(M,g_{\mathrm{ALG}})
$ equals the rank of the intersection matrix on $H^2(M,\del M)$.
\end{corollary}
\begin{proof} The intersection cohomology of $X$ can be calculated
using Mayer-Vietoris, so that the result follows from Theorem 1.
However, another approach may be more transparent. By
Theorem 3 and the signature formula (\ref{eq:sgnra}) 
it is enough to show that the fibration $\partial(M) \to B$ has 
$\tau$-invariant equal to $0$. But this follows from pp. 316-319 in
\cite{dai}, where it is shown that $\tau = 0$ on any fibration 
which admits a flat connection. This applies in the present situation
because over the one-dimensional base $B=S^1$ any connection is flat.
\end{proof}

In the examples of type $D_4$, the intersection matrix is the Cartan 
matrix of type $\hat{D}_4$,  \cite{cherkis-kapustin3}. Hence
in this case $L^2\calH^2(M,g_{\mathrm{ALG}})$ is four dimensional.

A parallel construction in \cite{cherkis-kapustin3} of certain
moduli spaces of solutions to Hitchin's equations (or equivalently 
Higgs bundles), yield  manifolds with hyperk\"ahler metrics 
$g_{\mathrm{Hit}}$ which have the same complex structure and underlying 
topology as the moduli spaces of periodic monopoles discussed above. 
A conjecture in \cite{cherkis-kapustin3} states that the corresponding
elements of these two classes of moduli spaces are in fact
isometric. For example, it is known that the moduli space of rank
$2$ parabolic Higgs bundles on $\CC P^1 \setminus \{p_1,p_2,p_3,p_4\}$ 
is an elliptic fibration (given by the Hitchin map) with one singular 
fiber of type $\hat{D}_4$. 

If this conjecture is valid in general, then Corollary 10 implies that 
for the $4$-dimensional moduli space of solutions to Hitchin's equations on a
cylinder, $L^2\calH^2(M,g_{\mathrm{Hit}}) \cong \mbox{\rm Im}
(H^2(M,\del M)\to H^2(M))$. This would be the first evidence, albeit 
indirect, for \cite[Conjecture 1]{hausel}.

\subsection{ALE toric hyperk\"ahler manifolds}
\label{toric}

Toric hyperk\"ahler manifolds have been defined and first studied
in \cite{bielawski-dancer}.  
An algebraic geometric account of the underlying varieties, with some
novel applications to combinatorics, is given in \cite{hausel-sturmfels}.

Let $\mbox{U}(1)^d$ act on $\HH^n$, preserving the hyperk\"ahler structure, 
and let $M_\xi=\HH^n////_{\xi}\mbox{U}(1)^d$ be a smooth toric hyperk\"ahler 
manifold of dimension $4n-4d$. The notation $X////_{\xi}G$ here denotes a
hyperk\"{a}hler quotient, see \cite{HKLR}.  This construction
determines a family of metrics on $M_{\xi}$ corresponding to the regular
values of the hyperk\"{a}hler moment map.  For any such value, consider
the family $M_{t\xi}$, $t > 0$. The asymptotics of the metrics in the family 
$M_{t\xi}$ are the same for $t \neq 0$ (i.e.\ these metrics
are quasi isometric, with increasing quasi-isometry constant as $t \to 0$).
As $t\to 0$, $M_{t\xi}$ degenerates to the singular space 
$M_0=\HH^n////_0 U(1)^d$. If we suppose that $M_0$ has only one isolated 
singularity, then the metrics in this family
maintain the same asymptotics at infinity even when $t=0$. In this case
$M_0$ is the cone over a $3$-Sasakian compact smooth manifold. This 
implies that $M_\xi$ is ALE. 

The question of when $M_0$ has only one isolated singularity is
intimately related to $3$-Sasakian geometry \cite{boyer-etal} and we 
quote a result from \cite[Theorem 4.1]{bielawski-dancer}:
$M_0$ has only one isolated singularity if and only if the action of
$\mbox{U}(1)^d$ on $\HH^n$ is unimodular (this means that the generic 
quotient $M_\xi$ is smooth) and generic (this means that the vector 
configuration described by the embedding $\mbox{U}(1)^d\subset \mbox{U}(1)^n$
is generic, see \cite{bielawski}). Now Theorem 1 and \cite{hausel-swartz}  
give
\begin{corollary} Suppose that the toric hyperk\"ahler manifold
$M_\xi$ is smooth and generic. Then 
\[
L^2\calH^{2n-2d}(M_\xi) \cong \mbox{\rm Im}\left(H^{2n-2d}(M_\xi,\del M_\xi)\to
H^{2n-2d}(M_\xi)\right) \cong H^{2n-2d}(M_\xi),
\]
and $L^2\calH^k(M_\xi) = 0$ in all other degrees. 
\end{corollary}
The fact that the Hodge cohomology is concentrated in the middle degree
is because $M_\xi$ has no cohomology above the middle dimension.
It is proven in \cite{hausel-swartz} that the intersection form on $H^{2n-2d}(M_\xi,\del M_\xi)$ is definite, 
which in the case of a smooth and 
generic toric hyperk\"ahler variety is consistent with Corollary 7. 
It follows that 
the forgetful map $ H^{2n-2d}(M_\xi,\del M_\xi)\to H^{2n-2d}(M_\xi)$ is
an isomorphism for any smooth toric hyperk\"ahler variety proving the last 
isomorphism in the above Corollary 11.

There are two extreme cases for a smooth generic toric hyperk\"ahler 
manifold $M_\xi$. One occurs when $d=n-1$, and these are just the ALE 
gravitational instantons of type $A_k$, which we have discussed earlier. 
The other extreme is when $d=1$, and then we obtain the Calabi metric
on $T^*\CC P^{n-1}$. From the argument above it has an ALE metric and its
Hodge cohomology is supported in the middle degree $2n-2$, where it is 
one-dimensional. An explicit generator for this space was found in
\cite{kim-etal}.

A closely related example is the ALE Ricci-flat K\"ahler metric on
$T^* S^n$, constructed by Stenzel in \cite{stenzel}.  Theorem 1 shows 
that there is a one-dimensional space of $L^2$ harmonic $n$-forms on 
that manifold when $n$ is even. For $n=2$ this is just the Eguchi-Hanson 
metric. For general $n=2k$, physicists have found explicit expressions for 
the $L^2$ harmonic $k$-form \cite{cvetic-etal}.

\subsection{$\mbox{Spin}(7)$ and $\mbox{G}_2$ metrics}
\label{special}
There has been recent interest amongst physicists to construct new
non-compact complete $\mbox{Spin}(7)$ and $\mbox{G}_2$ metrics, cf.\ 
\cite{gomis-etal}, and there have been predictions about the $L^2$ 
harmonic forms on such spaces. All known examples have fibred
boundary metrics, and so our results, Theorem 1, (\ref{sphereeven}) and 
(\ref{sphereodd}) can be used to check these predictions.
We mention just two examples. 

In fact, our Theorem 1 suggested that physicists look for an $L^2$ harmonic 
$3$-form on a particular example, an ALE $G_2$ metric on a rank $3$ real
vector bundle over $S^4$, constructed first in \cite{br-sa}.
We have as a simple corollary of Theorem 1:
\begin{corollary} The $\mbox{G}_2$ metric of \cite{br-sa}
on a rank $3$ real vector bundle over $S^4$ supports exactly a 
$1$-dimensional space of degree $3$ and a $1$-dimensional space of 
degree $4$ $L^2$ harmonic forms. 
\end{corollary}
Armed with the knowledge that such forms existed,
physicists \cite{gibbons-etal} were able to find their explicit
forms, see (2.18) of \cite{gibbons-etal} and also 
Footnote 4 in \cite{gibbons-etal}.

There is another example of a $\mbox{G}_2$ metric, constructed in 
\cite{gomis-etal}, which lives on $\RR^4\times S^3$.
It is ALF with $F=S^1$ and so our result (\ref{sphereodd}) implies that
\begin{corollary} There are no non-trivial $L^2$ harmonic forms on
the $G_2$ space of \cite{gomis-etal}.
\end{corollary}

A prediction coming from duality arguments between M-theory and type 
IIA string theory suggested the existence on this space of an $L^2$ 
harmonic $2$-form, or equivalently, a finite energy $\mbox{U}(1)$ Yang-Mills 
field, whose counterpart exists in the dual theory. This last corollary 
shows that this prediction fails; actually, already the methods of 
\cite{hitchin2} were used in \cite[Section 6]{gomis-etal} to establish the 
non-existence of $L^2$ harmonic $2$-forms on this $\mbox{G}_2$ manifold. 
Those authors call this the $\mbox{U}(1)$ puzzle.

\end{document}